\documentclass[a4paper,10pt,openright,twoside]{article}

\usepackage{setspace}
\usepackage{color}
\usepackage{latexsym}
\usepackage{amssymb,amsbsy,amsmath,amsfonts,amssymb,amscd}
\usepackage{epsfig, graphicx}
\usepackage{hyperref}
\usepackage{enumerate}
\usepackage{mathrsfs}
\usepackage{color}

\usepackage{float}
\usepackage{tikz}
\usepackage[all]{xy}
\usepackage{wrapfig}
\usepackage{enumitem}
\usepackage{multirow, array} 
\usepackage{setspace}
\usepackage{appendix}
\usepackage{fancyhdr}

\newcommand{\commentout}[1]{}
\newcommand{\R}{\mathbb{R}}

\newcommand {\Chi} {{\bf \raise 2pt \hbox{$\chi$}} }
\newcommand{\norme}[1]{\left\lVert#1\right\rVert}

\newcommand{\beq}{\begin{equation}}
\newcommand{\eeq}{\end{equation}}
\newcommand{\bea} {\begin{array}{rl}}
\newcommand{\eea} {\end{array}}
\newcommand{\bepa}{\left\{ \begin{array}{l}}
\newcommand{\eepa} {\end{array}\right.}
\newtheorem{theorem}{Theorem}[section]
\newtheorem{lemma}[theorem]{Lemma}
\newtheorem{definition}[theorem]{Definition}
\newtheorem{remark}[theorem]{Remark}
\newtheorem{proposition}[theorem]{Proposition}
\newtheorem{corollary}[theorem]{Corollary}
\providecommand{\abs}[1]{\lvert#1\rvert}

\numberwithin{equation}{section}

\newcommand{\qed}{{ \hfill
                      {\unskip\kern 6pt\penalty 500 \raise -2pt\hbox{\vrule\vbox to 6pt{\hrule width 6pt
                      \vfill\hrule}\vrule} \par}   }}

\title{\Large \bf Global-in-time classical solutions  and qualitative 
properties for the NNLIF neuron model with synaptic delay}
\author{Mar\' {\i}a J. C\'aceres
\thanks{Departamento de Matem\'atica Aplicada,
Universidad de Granada, 18071 Granada, Spain.
\texttt{caceresg@ugr.es}}
\and Pierre Roux\thanks{Laboratoire de biologie computationnelle et quantitative (LCQB), UMR 7238 CNRS, Universit\'e Pierre et Marie Curie, 75205 Paris Cedex 06, France et Laboratoire de Math\'{e}matiques d'Orsay (LMO), Universit\'{e} Paris-Sud, Paris-Saclay, Orsay, France.
	\texttt{pierre.roux@math.u-psud.fr}}
\and Delphine Salort \thanks{Laboratory of Computational and Quantitative Biology (LCQB),
	UMR 7238 CNRS, Universit\'e Pierre et Marie Curie, 75205 Paris Cedex 06, France.
	\texttt{delphine.salort@upcm.fr}}
\and Ricarda Schneider
\thanks{Departamento de Matem\'atica Aplicada,
Universidad de Granada, 18071 Granada, Spain.
\texttt{ricarda.schneider87@gmail.com}}
}
\date{\today}

\begin{document}
\maketitle           
\begin{abstract}
  The Nonlinear Noisy Leaky Integrate and Fire (NNLIF) model is widely
  used to describe the dynamics of neural networks after a diffusive
  approximation of the mean-field limit of a stochastic differential
  equation system.  When the total activity of the network has an
  instantaneous effect on the network, in the average-excitatory case,
  a blow-up phenomenon occurs.  This article is devoted to the
  theoretical study of the NNLIF model in the case where a delay in
  the effect of the total activity on the neurons is added.  We first
  prove global-in-time existence and uniqueness of classical
  solutions, independently of the sign of the connectivity parameter,
  that is, for both cases: excitatory and inhibitory. Secondly, we
  prove some qualitative properties of solutions: asymptotic
  convergence to the stationary state for weak interconnections and a
  non-existence result for periodic solutions if the connectivity
  parameter is large enough.  The proofs are mainly based on an
  appropriate change of variables to rewrite the NNLIF equation as a
  Stefan-like free boundary problem, constructions of universal
  super-solutions, the entropy dissipation method and Poinca\-r\'e's
  inequality.
\end{abstract}
\newline
\noindent
{\bf AMS Class. No}:  35K60, 82C31, 92B20




\section{Introduction}

Different scientific disciplines study the complex dynamics of neural
networks. Over the last decades, mathematicians have been particularly
interested in providing specific models to understand the beha\-vior
of neurons.
One popular approach is to tackle qualitative properties of networks  
via partial differential equations by deriving mean fields models from 
stochastic differential equations.
Depending on the choice of the intrinsic dynamics of  neurons and on
 the type of synaptic coupling one may obtain different 
models (see, among others, 
\cite{chevallier2015microscopic,chevallier2015mean,PPD2,Mischler_TENN_2018,
PeSa,PPD,RBW,delarue2015global,delarue2015particle,Caikinetic, NKKZRC10,CCP,
Henry:13,dumont2017mean, mischler2016kinetic}). 
In this article we assume that  the neurons  are described via their membrane
 potential and that when the membrane potential reaches a critical or \emph{threshold value} 
$V_F$, the neurons emit an action potential, also called spike, 
as a result of depolarization of their membrane, 
and their voltage values return to a \emph{reset value} $V_R$ ($V_R<V_F$). 
More precisely, we consider the following  PDE model 
(see  \cite{BrGe,brunel,BrHa} for its derivation):
 \begin{eqnarray}
\frac{\partial \rho}{\partial t}(v,t) +\frac{\partial}{\partial v}[(-v+\mu(t-D)\rho(v,t)]-a \frac{{\partial}^2 \rho}{\partial v^2}(v,t) =N(t)\delta(v-V_R), \ \ \ \ \  v\leq V_F, \label{delay_original}
\end{eqnarray}
where the function $\rho(\cdot,t)$ is the probability density of the
electric potential of a randomly chosen neuron at time $t$, $D\ge 0$
is the delay between emission and reception of the spikes (synaptic
delay), $V_R$ is the reset potential after firing, and
$V_F \in \mathbb{R }$ is the threshold potential.  The drift term
$\mu$ and the firing rate $N$ of the network are given by
\begin{eqnarray}\label{condiciones1}
\mu(t)=b_0+bN(t) \ \ \ \ \textrm{with} \ \ \ \ \ N(t)=-a \frac{\partial p}{\partial v}(V_F,t) \geq 0,
\end{eqnarray}
where $a>0$ is the diffusion coefficient, the parameter $b_0$ controls 
the strength of the external stimuli and can be either zero, 
positive or negative, and $b$ is  the  \emph{connectivity parameter}.
The neurons of the network can be either excitatory or inhibitory. 
This property is reflected in the NNLIF equation through the sign of $b$: 
$b>0$ for average-excitatory networks and $b<0$ for average-inhibitory ones.
\noindent The PDE \eqref{delay_original} is completed with initial and 
boundary conditions 
\begin{equation}\label{condiciones2}
N(t)=N^0(t)\geq 0, \quad \forall \ t\in[-D,0], \quad 
\rho(v,0)=\rho^0(v) \geq 0,  \quad \textrm{and} 
\quad \rho(V_F,t)=\rho(-\infty,t)=0.
\end{equation}
These boundary conditions imply that the following constraint is satisfied
$N^0(0) = -a\dfrac{\partial \rho^0}{\partial v}(V_F)$. 
Besides, for any classical solution $\rho$, the total mass is
conserved: if $\rho^0$ is a probability density, then that is also
true of $\rho$ at any positive time:
$
\int_{-\infty}^{V_F} \rho(v,t)\ dv= \int_{-\infty}^{V_F} \rho^0(v)\ dv = 1.$
The behaviour of the network, and, of course, of the solutions of the 
NNLIF equation depends strongly on the type of network considered, in terms of 
the sign of $b$.
In fact, in \cite{CCP} (among others) it was shown that there are some
situations in the case $D=0$, depending on the initial data and the
size of $b$, where the solutions for the average-excitatory case
cannot be global-in-time. The simulations therein suggest that the
blow-up phenomenon is reflected in a divergence in finite time of the
firing rate.  Following these observations, in \cite{CGGS} a criterion
for the maximal time of existence for classical solutions was
derived. Essentially, it ensures that the solutions exist while the firing
rate is finite. Besides, it was obtained that classical solutions for
the average-inhibitory case are always global-in-time. In
\cite{carrillo2014qualitative}, some qualitative properties were
proved in the case $D=0$: uniforms bounds in $L^\infty$ for the
average-inhibitory case and asymptotic convergence toward steady
states for small connectivity parameters.

Nevertheless, all of these previous works analyze a version of the
NNLIF equation with no synaptic delay. The synaptic delay $D$ is the
short period of time that passes since a nerve impulse is sent from a
presynaptic neuron until it finally reaches the postsynaptic
neuron. The numerics done in \cite{CS17-2} suggest that the synaptic
delay in the NNLIF equation prevents the blow-up of the firing
rate. Instead, the firing rate is seen to converge to a stationary
state, oscillate or increase. In fact, at the microscopic level, it
has already been proved that the solutions are global-in-time for the
delayed NNLIF model \cite{delarue2015particle}.  Moreover, for the
population density model of IF neuron with jumps, which arises from
the same microscopic approximation as the NNLIF model, it was shown in
\cite{Henry:13} that the firing rate blows up in finite time in some
situations, but in \cite{Henry:12} it was proved that this blow-up
disappears if a synaptic delay is considered.

In the present paper we deal with the delayed NNLIF equation, which is
a modification of the NNLIF model presented in \cite{CCP} and
\cite{CGGS} at the level of the drift term, which includes a synaptic
delay $D>0$.  We prove the global-in-time existence of classical
solutions for both the average-inhibitory and the average-excitatory
cases.  Moreover, we analyse the long time behaviour of these
solutions for a small connectivity parameter, and we provide a
non-existence result of periodic-in-time solutions for a large enough
connectivity parameter.

For simplicity, sometimes we will suppose that $a=1$ and
$V_F=0$. Nevertheless, these hypotheses are not really a constraint,
since we can transform the general equation into one satisfying the
restriction by defining a new density $\bar \rho$ as follows:
\begin{equation}\label{traslado}
\bar \rho(v,t)=\sqrt{a}\rho(\sqrt{a} v+ V_F,v).
\end{equation}
Also for simplicity, sometimes we will assume that $b_0=0$, but again
we can do it without loss of generality, since we can pass from the
general equation for $b_0 \ne 0$ to an equation with $b_0=0$ by
translating the voltage variable $v$ by the factor $b_0$.

\

The structure of the paper is as follows: the second and third sections are 
devoted to the study of the Cauchy problem of Equation \eqref{delay_original}. 
Taking into account  the presence of the delay, we adapt
the strategy of \cite{CGGS} in the second section, and rewrite the 
system \eqref{delay_original}-\eqref{condiciones1}-\eqref{condiciones2} as  
a free boundary Stefan-like problem with a nonstandard right hand side 
consisting of a Dirac Delta source term. 
We also provide a definition of classical solutions for the new
problem, and give some a priori properties for them.
The third section is devoted to the proof of local
existence and uniqueness of classical solutions for the
Stefan-like problem. This is done through a fixed point argument for
an integral formulation of the Stefan-like equation. Besides, it is
shown how this result can be extended to the Fokker-Planck Equation
\eqref{delay_original}.  The fourth section contains one of the main
results. There we prove the global existence and uniqueness of
classical solutions for both the average-inhibitory and the
average-excitatory case when $D>0$. First, we extend to the case at
hand the characterization of the maximal time of existence of the
solutions in terms of the size of the firing rate, provided in
\cite{CGGS} for the case without delay ($D=0$). Mainly, this ensures that
local solutions exist and are unique as long as the firing rate is
finite. As in the case $D=0$, the global existence for the
average-inhibitory case is derived, showing that every solution
defined until a certain time $t_0$ can be extended up to a short (but
uniform) time $\varepsilon$, since the firing rate up to this
additional time $t_0+\varepsilon$ is uniformly bounded.  However, for
the average-excitatory case, this uniform bound for the firing
rate is not arrived at easily, even with delay (see \cite{CS17-2} for
numerical results).  We overcome this difficulty with a super-solution
for the delayed NNLIF equation. We do not prove a uniform bound of the
firing rate, but we show that for every maximal time of existence the
firing rate of the solution is bounded. Therefore, the criterion of
maximal time of existence gives us a contraction and consequently we
prove global existence of the solution.  The fifth section is devoted
to studying the long time behavior of the system
\eqref{delay_original}-\eqref{condiciones1}-\eqref{condiciones2}. We
show exponentially fast convergence to the steady state of the
solutions if the connectivity parameter $b$ is small, extending the
results of \cite{CCP} and \cite{carrillo2014qualitative}. This result
is achieved, for both the average-excitatory and the
average-inhibitory case, by means of the entropy method, a
Poincar\'e's inequality and suitable $L^2$ estimates of the firing
rate. We also study restrictions on the existence of time-periodic
solutions, although it is still an open question whether they exist
or not in other situations.  Finally, the appendix contains technical
tools used extensively in the paper: the Poincar\'e-like inequality
and the entropy equality.


\section{The equivalent free boundary Stefan problem}

In this section we rewrite Equation (\ref{delay_original}) 
as a free boundary Stefan problem with a nonstandard right hand side 
as in \cite{CGGS} was done without delay.  
With that purpose we  perform the three changes of variables presented below.
Afterwards, we write the final expression of the equivalent equation,
 we introduce the notion of classical solution for it and remember 
some basic a priori properties for this kind of solutions.
In all this section, we assume $D>0$, $a=1$ and $V_F=0$.
\begin{enumerate}
\item {\bf A first change of variables.} We introduce the following change of variables, which  has been widely studied in \cite{carrillo2000asymptotic}:
\begin{eqnarray}\label{cambio1}
y=e^t v, \quad \tau=\frac{1}{2}(e^{2t}-1).
\end{eqnarray}
Therefore, denoting by $\alpha(\tau)=(\sqrt{2\tau+1})^{-1}$,
\begin{align}\label{equivalencia}
 t=-\log(\alpha(\tau)), \quad v=y\alpha(\tau),
\end{align}
and we define
\begin{eqnarray}
w(y, \tau)=\alpha(\tau) \rho \big(y \alpha(\tau), -\log(\alpha(\tau))\big). \label{w_nuevo}
\end{eqnarray}
 Differentiating $w$ with respect to $\tau$, and using that $\rho$ 
is a solution of \eqref{delay_original}, yields
\begin{eqnarray}
\dfrac{\partial w}{\partial \tau}(y,\tau)
 &=& \alpha'(\tau)\rho(y\alpha(\tau), -\log(\alpha(\tau)))\nonumber\\
& \ &+y\alpha'(\tau)\alpha(\tau)\dfrac{\partial \rho}{\partial v}(y\alpha(\tau),-\log(\alpha(\tau)))-\alpha'(\tau)\rho_t(y\alpha(\tau),-\log(\alpha(\tau)))\nonumber\\
 & = & 
-\alpha'(\tau)\dfrac{\partial^2 \rho}{\partial v^2}(y\alpha(\tau),-\log(\alpha(\tau))) \\ \label{w_t2}
&  & +\alpha'(\tau)\mu(-\log(\alpha(\tau))-D)\dfrac{\partial \rho}{\partial v}(y\alpha(\tau), -\log(\alpha(\tau))) \nonumber
\\ 
& & -\alpha'(\tau)N(-log(\alpha(\tau)))\delta(y\alpha(\tau)-V_R). \nonumber
\end{eqnarray}
Finally, taking into account that $-\alpha'(\tau)=\alpha^3(\tau)$ and 
\[	\dfrac{\partial w}{\partial y}(y,\tau)=\alpha^2(\tau)\dfrac{\partial \rho}{\partial v}(y\alpha(\tau),-\log(\alpha(\tau))),\]
\[	\dfrac{\partial^2 w}{\partial y^2}(y,\tau)=\alpha^3(\tau)\dfrac{\partial^2 \rho}{\partial v^2}(y\alpha(\tau),-\log(\alpha(\tau))), 
\]
we obtain
\begin{eqnarray}
\dfrac{\partial w}{\partial \tau}(y,\tau)=\dfrac{\partial^2 w}{\partial y^2}(y, \tau)-\alpha(\tau)\mu(t-D)\dfrac{\partial w}{\partial y}(y,\tau)+M(\tau)\delta\left(y-\frac{V_R}{\alpha(\tau)}\right), \label{eq_w}
\end{eqnarray}
where $M(\tau)=-\dfrac{\partial w}{\partial y}(0,\tau)=\alpha^2(\tau)N(t)$, and we use that  $\mu(t-D)=\mu(-\log(\alpha(\tau))-D)$ due to \eqref{equivalencia}. 
Therefore, $w$ satisfies the equation
\[ 
\left\{ \begin{array}{ll}
\dfrac{\partial w}{\partial \tau}(y,\tau) = \dfrac{\partial^2 w}{\partial y^2}(y,\tau) - \alpha(\tau)\mu(t-D)\dfrac{\partial w}{\partial y}(y,\tau) +M(\tau)\delta_{\frac{V_R}{\alpha(\tau)}}( y) & \qquad y\in]-\infty,0], \quad \tau\in\R_+, \\
M(\tau)=-\dfrac{\partial w}{\partial y}(0,\tau)& \qquad \tau\in\R_+,  \\
w(-\infty,\tau)=w(0,\tau)=0 & \qquad \tau\in\R_+,\\
\displaystyle w(y,0)=w_I(y)  &\qquad  y\in(-\infty,0].
\end{array}\right.
\]
\item {\bf A second change of variables.}
With the change of variable (denoting $t_\omega= -\log(\alpha(\omega))$), 
\begin{equation}
x =y - \int_{0}^{\tau}\mu\left(t_\omega-D\right)\alpha(\omega)d \omega
=
 y - \int_{0}^{\tau}\mu\left(-\log(\alpha(\omega))-D\right)\alpha(\omega)
d \omega, 
\label{cambio2}\end{equation} 
the function $u$ defined by $u(x,\tau)=w(y,\tau)$ satisfies
\begin{equation}
\left\{\begin{array}{ll}
\displaystyle \dfrac{\partial u}{\partial \tau}(x,\tau) = \dfrac{\partial^2 u}{\partial x^2}(x,\tau) +M(\tau)\delta_{s(\tau)+\frac{V_R}{\alpha(\tau)}}( x) &\qquad x\in]-\infty,s(\tau)], \tau\in\R_+, \\
\displaystyle s(\tau) = - \int_{0}^{\tau}\mu\left(-\log(\alpha(\omega))-D\right)\alpha(\omega)d \omega & \qquad \tau\in\R_+,\\
\displaystyle M(\tau)=-\dfrac{\partial u}{\partial x}(s(\tau),\tau)& \qquad \tau\in\R_+,\\
\displaystyle u(-\infty,\tau)=u(s(\tau),\tau)=0& \qquad \tau\in\R_+,\\
\displaystyle u(x,0)=u_I(x)  &\qquad  x\in(-\infty,0].
\end{array}\right.  \label{EqnEq}
\end{equation}

\item {\bf A third change of variables.}
In the system  \eqref{EqnEq} the boundary $s(\tau)$ depends on
$\mu\left(t-D\right)=b_0+bN\left(t-D\right)$, and it, in turn, on $M$. 
Therefore, we have to remove the explicit  $t$ dependency in the term $N(t-D)$.
In the case $D=0$ (see \cite{CGGS}) it is easy since 
$\mu\left(t-D\right)=b_0+bN\left(t-D\right)=b_0+bM(\tau)\alpha^{-2}(\tau)$. 
However, if $D>0$ the relation is more involved, because
$t=\frac{1}{2}\log(2\tau+1)$ and $\tau=\frac{1}{2}\left(e^{2t}-1\right)$, 
and thus, if we consider the 
time $t-D$, there is a related $\tau_D$:
$\tau_D=\frac{1}{2}(e^{2(t-D)}-1)$, which makes complicated the handling of $D$.
We overcome this difficulty considering 
 $\bar D=(1-e^{-2D})>0$. With this choice $\tau_D=\tau\left(1-\bar{D}\right)-
\frac{\bar D}{2}$ and we obtain:
\begin{eqnarray}
N(t-D)&=&\alpha^{-2}(\tau_D)M(\tau_D)=
\alpha^{-2}\left((1-\bar D)\tau-\frac{1}{2}\bar D\right)
M\left((1-\bar D)\tau-\frac{1}{2}\bar D\right), \label{N_delay}
\end{eqnarray}
which depends only on $\tau$.
Thus, the initial synaptic delay $D$ is translated into the delay 
$\bar D$, which is scaled between $0$ and $1$, being $\bar D=0$ if $D=0$ 
and $\bar D=1$ if $D=\infty$.
In this way, using (\ref{N_delay}) we can rewrite $s(\tau)$ in terms of $M(\tau)$,  avoiding the dependence in $t$ 
\begin{eqnarray}
s(\tau)&=&-b_0(\sqrt{2\tau+1}-1)-b\int_0^\tau N(t_\omega-D)\alpha(\omega) \ d\omega  \label{s_sin_cv}\\
& = &-b_0(\sqrt{2\tau+1}-1)-
b\int_0^\tau \alpha^{-2}\left((1-\bar D)\omega-\frac{1}{2}\bar D\right)
M\left((1-\bar D)\omega-\frac{1}{2}\bar D\right)\alpha(\omega) \ d\omega. \nonumber
\end{eqnarray}
The change of variables $z=(1-\bar D)\omega - \frac{1}{2}\bar D$ yields
\begin{eqnarray} \label{cambio_s}
s(\tau)&=&-b_0(\sqrt{2\tau+1}-1)-\frac{b}{\sqrt{1-\bar D}}
\int_{-\frac{1}{2}\bar D}^{(1-\bar D)\tau-\frac{1}{2}\bar D} M(z)\alpha^{-1}(z)  
\ dz. \label{s_con_cv}
\end{eqnarray}
Taking into account that $M^0(\tau)=\alpha^2(\tau)N^0\left(\frac{1}{2}\log(2\tau+1)\right)$,  finally leads to the following equi\-va\-lent Stefan-like equation
\begin{eqnarray}
\left \{
\begin{array}{ll}
\dfrac{\partial u}{\partial \tau}(x,\tau) =  \dfrac{\partial^2 u}{\partial x^2}(x, \tau)+M(\tau)\delta_{s(\tau)+\frac{V_R}{\alpha(\tau)}}(x), & x<s(\tau), \tau>0,   \\ 
s(\tau)=-b_0(\sqrt{2\tau+1}-1)-\frac{b}{\sqrt{1-\bar D}}
\displaystyle \int_{-\frac{1}{2}\bar D}^{(1-\bar D)\tau-\frac{1}{2}\bar D} 
M(z)\alpha^{-1}(z) \ \ dz,&   \tau>0, \\  
M(\tau)=-\dfrac{\partial u}{\partial x}(s(\tau),\tau),&   \tau>0, \\ 
M(\tau)=M^0(\tau)>0,&    \tau \in (-\frac{\bar D}{2},0], \\  
u(-\infty, \tau)=u(s(\tau),\tau)=0, &  \tau>0, \\  
u(x,0)=u_I(x),  &  x<0.  
\label{eq_u_definitiva}
\end{array}
\right.
\end{eqnarray}
where  $\bar D \in [0,1)$ and $\alpha(\tau)=\frac{1}{\sqrt{2\tau+1}}$. 
Let us remark that this problem is well defined since 
$\alpha(\tau)\in \R^+, \ \forall \ \tau > -\frac{1}{2}$.
Also, we note that if $\bar D=0$, the system \eqref{eq_u_definitiva} 
reduces to the system studied in \cite{CGGS}. 
\end{enumerate}

We conclude this section with the notion of classical solutions 
for this kind of system and with some a priori properties, 
that will be useful for the rest  of the computations of the present work.
\begin{definition}[Classical solutions for the Stefan-like problem]\label{def_u}
	Let $u^0(x)$ be a non-negative $C^0((-\infty,0])\cap C^1((-\infty, V_R)\cup(V_R,0]) \cap L^1((-\infty,0))$ function such that $u^0(0)=0$. Suppose that  $\frac{d u^0}{dx}$ vanish at $-\infty$ and admits finite left and right limits at $V_R$. We say that $u$ is a classical solution of \eqref{eq_u_definitiva} 
(equivalently \eqref{EqnEq}) with initial datum $u^0$ on the interval $J=[0,T)$ or $J=[0,T]$, for a given $T>0$ if:
	\begin{enumerate}
		\item $M(\tau)$ is a continuous function for all $\tau \ \in \ J$,
		\item $u$ is continuous in the region $\{(x,\tau): -\infty<x\leq s(\tau), \tau \in J\}$ and for all $\tau\in J$, $u\in L^1((-\infty,s(\tau)))$,
		\item $\partial_{xx}u$ and $\partial_\tau u$ are continuous in the region $\{(x,\tau): -\infty<x<s_1(\tau), t \in J \setminus \{0\}\} \ \cup \  \{ (x,\tau): s_1(\tau)<x<s(\tau), \tau \in J\setminus\{0\}\}$,
		\item If we denote  $s_1(\tau):=s(\tau)+\dfrac{V_R}{\alpha(\tau)}$, then $\partial_xu(s_1(\tau)^-, \tau)$, $\partial_xu(s_1(\tau)^+,\tau)$, $\partial_xu(s(\tau)^-, t)$ are well defined,
		\item $\partial_xu$ vanishes at $-\infty$,
		\item Equations \eqref{eq_u_definitiva} (equivalently \eqref{EqnEq}) are satisfied.
	\end{enumerate}
\end{definition}
Obviously  we observe that this definition includes the notion of solution
given in \cite{CGGS} for the case without transmission delay ($D=\bar D=0$).

We gather in the following lemma some a priori properties
of these solutions.
\begin{lemma}[A priori properties]
	Let $u$ be a solution to  \eqref{eq_u_definitiva} (equivalently \eqref{EqnEq}) in the sense of the previous definition. Then:
	\begin{enumerate}
		\item The mass is conserved: $\int_{-\infty}^{s(\tau)} u(x,\tau)dx=\int_{-\infty}^{0} u^0(x)dx$, $\forall \  t>0$.
		
		\item The flux across the moving point $s_1(\tau):=s(\tau)+\dfrac{V_R}{\alpha(\tau)}$ is exactly the strength of the source term:
		\[
			M(\tau):=-\partial_xu(s(\tau),\tau)=\partial_xu(s_1(\tau)^-,\tau)-\partial_xu(s_1(\tau)^+,\tau).
		\]
		\item If $b_0\leq0$ and $b<0$ (respectively, $b_0\ge0$ and $b>0$), the free boundary $s(\tau)$ is a monotone increasing (respectively, decreasing) function of time.
	\end{enumerate}
\end{lemma}
{\bf Proof.} The proof of properties $1.$ and $2.$ is exactly the same as in \cite{CGGS}[Lemma 2.3], because it does not take into account the expression of $s(t)$. Property $3.$ is obvious with the form of $s$.


\section{Local existence and uniqueness}
\label{sec: local}

In this section we introduce an implicit integral equation for $M$. 
Then, thanks to the form of that equation, it  will be possible to solve
it for local time  using a fixed point argument. 
Besides, since we will be able to prove that the fixed point function is 
a contraction, we will also get the local uniqueness of $M$. 
This is useful, since once $M$ is known,   \eqref{eq_u_definitiva} 
(equivalently \eqref{EqnEq}) decouples, and $u$ can be calculated easily 
by a Duhamel's formula.

The inclusion of the transmission delay, $D>0$, in the model produces 
that the NNLIF equation \eqref{delay_original}
becomes linear  for $t<D$, since the drift and the diffusion terms
depend on the initial condition, instead of the firing rate $N(t)$.
This fact is translated to the   System \eqref{eq_u_definitiva} 
(equivalently \eqref{EqnEq})  when $\tau \leq \frac{1}{2}(e^{2D}-1)=
\frac{\bar D}{2(1-\bar D)}$, which means $t \leq D$ 
in the original time variable. 
Thus, the boundary is free with a delay but it is constrained on 
any short period of time, since on every time interval of size less 
than $\frac{1}{2}(e^{2D}-1)=\frac{\bar D}{2(1-\bar D)}$, systems 
\eqref{EqnEq} and \eqref{eq_u_definitiva} are equivalent to a linear system
 with non-standard right-hand side. To be more precise,
 we can solve them by studying the following equivalent
linear system:
\begin{equation}
\left\{\begin{array}{ll}
\displaystyle \dfrac{\partial u}{\partial \tau}(x,\tau) = \dfrac{\partial^2 u}{\partial x^2}(x,\tau) +M(\tau)\delta_{s(\tau)+\frac{V_R}{\alpha(\tau)}}( x) &\qquad x\in]-\infty,s(\tau)], \tau\in\R_+, \\
\displaystyle s(\tau) = - \int_{0}^{\tau}I(\omega)d \omega & \qquad \tau\in\R_+,\\
\displaystyle u(-\infty,\tau)=u(s(\tau),\tau)=0& \qquad \tau\in\R_+,\\
\displaystyle M(\tau)=-\dfrac{\partial u}{\partial x}(s(\tau),\tau)& \qquad \tau\in\R_+ , \\  
u(x,0)=u_I(x),  & \qquad x<0.  

\end{array}\right.  \label{AbstractEq}
\end{equation}
where $I\in\mathscr{C}^0([0,+\infty))$ is an abstract input function.
In our case, 
 the system \eqref{EqnEq} (equivalently \eqref{eq_u_definitiva})
can be written as System  \eqref{AbstractEq} with
$I(\omega)=\mu\left(-\log(\alpha(\omega))-D\right)\alpha(\omega)$,
 on every time interval of size less 
than $\frac{1}{2}(e^{2D}-1)=\frac{\bar D}{2(1-\bar D)}$,
since $I$ does not depend on $u$ on this time, but on the values
on the previous time intervals. Therefore,  
System \eqref{EqnEq} (equivalently \eqref{eq_u_definitiva})
can be considered as a linear system on every time interval of size less 
than $\frac{1}{2}(e^{2D}-1)=\frac{\bar D}{2(1-\bar D)}$.

The notion of solution of definition \ref{def_u} and the 
\textit{a priori} properties still apply for equation \eqref{AbstractEq} if 
we assume $I$ such that the free bondary $s(\tau)$ is a  monotone function
on time.

Denoting $G(x,\tau,\xi,\eta) = \dfrac{1}{\sqrt{4\pi (\tau-\eta) }} e^{-\dfrac{(x-\xi)^2}{4(\tau-\eta)}}$, we prove the following theorems.
\begin{theorem}
	If the function $u$ is solution of \eqref{AbstractEq} in the sense of definition \eqref{def_u}, then the continuous function $M$ satisfies
	\begin{multline}  M(\tau) = - 2 \int_{-\infty}^{0} G(s(\tau),\tau,\xi,0)\dfrac{d u^0}{d x}(\xi)d \xi \\ + 2\int_{0}^{\tau} M(\eta)\dfrac{\partial}{\partial x} G(s(\tau),\tau,s(\eta),\eta)d \eta \\  - 2\int_{0}^{\tau} M(\eta)\dfrac{\partial}{\partial x} G(s(\tau),\tau,s(\eta)+\dfrac{V_R}{\alpha(\eta)},\eta)d \eta \label{ConditionM} \end{multline}
\end{theorem}
{\bf Proof.}
	The proof is exacly the same as in \cite{CGGS}. Therefore, we only sketch it. The main idea is to use the following Green's identity:
$ \dfrac{\partial}{\partial \xi} \left( G\dfrac{\partial u}{\partial \xi} - u\dfrac{\partial G}{\partial \xi} \right) - \dfrac{\partial }{\partial \tau}\big( Gu \big)  = 0$. 
	The result is obtained by integration of this identity in the regions $\{ \xi\in(-\infty, s(\eta)+V_R{\alpha(\eta)}^{-1}),\eta\in(0,\tau)\}$ and $\{ \xi\in(s(\tau)+V_R{\alpha(\tau)}^{-1},s(\tau)), \ \eta\in(0,\tau)\}$ and addition of the results. Computations are then done on the different terms to express them in function of $M$, $G$, $\partial_x G$ and $u^0$ under integrals.
\qed

\begin{theorem}
	Let $u^0(x)$ be a non-negative $C^0((-\infty,0])\cap 
        C^1\left((-\infty, V_R)\cup(V_R,0]\right) \cap L^1((-\infty,0))$ 
        function such that $u^0(0)=0$, and suppose that $\dfrac{d u^0}{dx}$ 
        vanishes at $-\infty$ and admits finite left and right limits at $V_R$.
        Then, there exists $T\in\R_+$ and a unique function 
        $M\in \mathscr{C}^0([0,T])$ that satisfies 
        \eqref{ConditionM} on $[0,T]$.
\end{theorem}

{\bf Proof.}
	Let $\sigma \in\R_+$, 
        $m:=1+2 \sup_{x \in (-\infty,V_R)~\cup~(V_R,0]} 
        \left|\dfrac{d u^0}{d x}(x)\right|$ and consider the space
	\[  C_{\sigma,m} = \{  M\in \mathscr{C}^0([0,\sigma]) \ | 
        \ \norme{M}_{\infty}:= \sup_{\tau\in [0,\sigma]}|M(\tau)|
        \leqslant m \}. \]
	We define in this space the functional
	\begin{eqnarray}  \mathcal{T}(M)(\tau)&:=&  
        - 2 \int_{-\infty}^{0} G(s(\tau),\tau,\xi,0)\dfrac{d u^0}{d x}(\xi)d 
        \xi  + 2\int_{0}^{\tau} M(\eta)\dfrac{\partial G}{\partial x} 
        (s(\tau),\tau,s(\eta),\eta)d \eta 
       \nonumber
       \\ & - & 2\int_{0}^{\tau} M(\eta)
       \dfrac{\partial G}{\partial x} (s(\tau),\tau,s(\eta)+
       \dfrac{V_R}{\alpha(\eta)},\eta)d \eta.
       \label{FonctionnelleM} 
       \end{eqnarray}
       The proof of the theorem is obtained if we show that this functional 
       has a unique fixed point. To do so, we start proving that for 
       $\sigma$ small enough $\mathcal{T}: C_{\sigma,m} \to C_{\sigma,m}$.\\
	On the one hand,
	\[  \begin{array}{rcl}
	\displaystyle \norme{ 2 \int_{-\infty}^{0} G(s(\cdot),\cdot,\xi,0)\dfrac{d u^0}{d x}(\xi)d \xi  }_\infty&=&\displaystyle 2 \sup_{x \in ]-\infty,V_R[~\cup~]V_R,0]} \left|\dfrac{d u^0}{d x}(x)\right| \norme{ \int_{-\infty}^{0} G(s(\cdot),\cdot,\xi,0)d \xi }_\infty\\
	&=&\displaystyle  2 \sup_{x \in ]-\infty,V_R[~\cup~]V_R,0]} \left|\dfrac{d u^0}{d x}(x)\right|.
	\end{array}  \]
	On the other hand, the positive valued applications
	\[  \phi_1 : \tau \mapsto 2m\int_{0}^{\tau} \left| \dfrac{\partial G}{\partial x} (s(\tau),\tau,s(\eta),\eta) \right|d \eta,\]
	and
	\[ \phi_2 : \tau \mapsto 2m\int_{0}^{\tau} \left|\dfrac{\partial G}{\partial x} (s(\tau),\tau,s(\eta)+\dfrac{V_R}{\alpha(\eta)},\eta)\right|d \eta, \]
	are continuous on $(0,\sigma]$. A direct computation gives
	\[ \left| \dfrac{\partial G}{\partial x} (s(\tau),\tau,s(\eta),\eta)\right| = \dfrac{1}{2\sqrt{4\pi}}\dfrac{|s(\tau) - s(\eta)|}{|\tau-\eta|^{\frac{3}{2}}}e^{-\dfrac{\big(s(\tau)-s(\eta)\big)^2}{4(\tau-\eta)}}  \]
	Since the function $I$ is bounded on every compact set 
      (because it is continuous), we can denote 
    $I_0=\sup_{\omega\in[0,\sigma]}|I(\omega)|$, and we have
	$|s(\tau)-s(\eta)| = 
        \left|\int_{\eta}^{\tau}I(\omega)d \omega \right| \leqslant 
      I_0 \left| \tau - \eta \right|$. 
	Substituting in the previous expression and putting back the integral, we get
	$\displaystyle \phi_1(\tau) \leqslant \displaystyle 
        \dfrac{mI_0}{\sqrt{4\pi}}\int_{0}^{\tau} 
         \dfrac{1}{\sqrt{\tau - \eta}} d \eta
	 = \dfrac{mI_0}{\sqrt{4\pi}}\sqrt{\tau}$.
	We directly get $\lim_{\tau\to0}\phi_1(\tau) = 0$ and $\phi_1$ is continuous on $[0,\sigma]$, with $\phi_1(0)=0$.
	Similarly, we have
	\[  |s(\tau)-s(\eta)-\dfrac{V_R}{\alpha(\tau)}| \geqslant |V_R| - I_0 | \tau - \eta |\geqslant |V_R| - I_0\tau,  \]
	and using the inequality $ ze^{-z^2} \leqslant e^{-z^2/2}$, 
we can write\\
	\[
	\displaystyle \phi_2(\tau) = \displaystyle  \dfrac{1}{2\sqrt{4\pi}} \int_{0}^\tau \dfrac{|s(\tau) - s(\eta) - \frac{V_R}{\alpha(\eta)}|}{|\tau-\eta|^{\frac{3}{2}}}e^{-\dfrac{\big(s(\tau)-s(\eta) - \frac{V_R}{\alpha(\tau)}\big)^2}{4(\tau-\eta)}}  d \eta, \]
	and then
	\[ \begin{array}{rcl}
	\displaystyle \phi_2(\tau)&\leqslant& \displaystyle \dfrac{1}{\sqrt{4\pi}} \int_{0}^\tau \dfrac{1}{\tau-\eta} e^{-\dfrac{\big(s(\tau)-s(\eta) - \frac{V_R}{\alpha(\tau)}\big)^2}{8(\tau-\eta)}} d \eta 
	\leqslant \displaystyle \dfrac{1}{\sqrt{4\pi}} \int_{0}^\tau \dfrac{1}{\tau-\eta} e^{-\dfrac{\big( |V_R| - I_0\tau \big)^2}{8(\tau-\eta)}} d \eta.\end{array}  \]
	When $\tau$ is small enough, $|V_R| - I_0\tau \geqslant \frac{1}{2}|V_R|$ and we have, making the change of variable $z=\frac{|V_R|}{2\sqrt{8(\tau-\eta)}}$,
	\[ \begin{array}{rcl}
	\displaystyle \phi_2(\tau) \leqslant \displaystyle \dfrac{1}{\sqrt{4\pi}} \int_{0}^\tau \dfrac{1}{\tau-\eta} e^{-\dfrac{\big( |V_R| - I_0\tau \big)^2}{8(\tau-\eta)}} d \eta
	&\leqslant & \displaystyle \dfrac{1}{\sqrt{4\pi}} \int_{0}^\tau \dfrac{1}{\tau-\eta} e^{-\dfrac{V_R^2}{32(\tau-\eta)}} d \eta\\
	&=& \displaystyle \dfrac{1}{\sqrt{\pi}}\int_{\frac{|V_R|}{2\sqrt{8\tau}}}^{+\infty} \dfrac{1}{z} e^{-z^2}d z.\end{array}  \]
	we also get $\lim_{\tau\to0}\phi_2(\tau) = 0$ and $\phi_2$ is continuous on $[0,\sigma]$, with $\phi_2(0)=0$.\\ 
	Thus, there exists an interval $[0,\zeta]\subset[0,\sigma]$ such that
	\[  \sup_{\tau\in[0,\zeta]} |\phi_1(\tau)| < \dfrac{1}{2} \qquad \mathrm{and} \qquad  \sup_{\tau\in[0,\zeta]} |\phi_2(\tau)| < \dfrac{1}{2}.  \]
	Hence, we have, for this $\zeta\in\R_+$ and for all $M\in C_{\zeta,m}$,
	\[  \norme{2\int_{0}^{\cdot} M(\eta)\dfrac{\partial G}{\partial x} (s(\cdot),\cdot,s(\eta),\eta)d \eta}_\infty \leqslant  \sup_{\tau\in[0,\zeta]}2m\int_{0}^{\tau} \left| \dfrac{\partial G}{\partial x} (s(\tau),\tau,s(\eta),\eta) \right|d \eta < \dfrac{1}{2}, \]
	and also
	\[ \norme{2\int_{0}^{\cdot} M(\eta)\dfrac{\partial G}{\partial x} (s(\cdot),\cdot,s(\eta)+\dfrac{V_R}{\alpha(\eta)},\eta)d \eta}_\infty < \dfrac{1}{2}.  \]
	Collecting the previous bounds, we obtain, for $\sigma$ small enough and depending only on $m$ and $I$, thus depending only on $\dfrac{d u^0}{d x}$ and $g$, that the map $\mathcal{T}$ is defined as a function from $C_{\sigma,m}$ into $C_{\sigma,m}$.

Finally, we prove that for this choice of $\sigma$ the
functional $\mathcal{T}$ is a contraction. We have for all $M,N\in C_{\sigma,m}$, for all $\tau\in[0,\sigma]$,
	\begin{multline} |\mathcal{T}(M)(\tau)-\mathcal{T}(N)(\tau)| \leqslant \dfrac{1}{m}\norme{M-N}_\infty\bigg(2m\int_{0}^{\tau}\left|\dfrac{\partial}{\partial x} G(s(\tau),\tau,s(\eta),\eta) \right|d \eta \\ + 2m\int_{0}^{\tau}\left|\dfrac{\partial}{\partial x} G(s(\tau),\tau,s(\eta)+\dfrac{V_R}{\alpha(\eta)},\eta)d \eta \right| \bigg) \\ \leqslant \dfrac{1}{m}\Big(\sup_{\tau\in[0,\sigma]} |\phi_1(\tau)|  + \sup_{\tau\in[0,\sigma]} |\phi_2(\tau)|\Big)\norme{M-N}_\infty.  \end{multline}
	As $m\geqslant1$, we have
	\[  \dfrac{1}{m} \Big( \sup_{\tau\in[0,\sigma]} |\phi_1(\tau)|  + 2\sup_{\tau\in[0,\sigma]} |\phi_2(\tau)| \Big) \\ < \dfrac{1}{2m} + \dfrac{1}{2m} = \dfrac{1}{m} \leqslant 1. \]
	Thus $\mathcal{T}$ is a contraction on the complete metric space $C_{\sigma,m}$. It admits a unique fixed point $M$ in that space.
\qed
\begin{theorem}[Local existence of the linear problem]\label{SolMax_u_linear}
	Let $u^0(x)$ be a non-negative  function in 
        $C^0((-\infty,0])\newline \noindent
       \cap  C^1\left((-\infty, V_R)\cup(V_R,0]\right) 
        \cap L^1((-\infty,0))$, such that $u^0(0)=0$ and suppose 
        that  $\frac{d u^0}{dx}$ vanishes at $-\infty$ and admits finite 
        left and right limits at $V_R$. Then, there exists a unique maximal 
        classical solution $u$ for  the problem \eqref{AbstractEq}.
\end{theorem}
{\bf Proof.}
The proof  is omitted, since it is performed as in \cite{CGGS} [Corollary 3.3].
Let us only point out that, once $M$ is known the equation for $u$ decouples,
 and $u$ can be calculated via the Duhamel's formula
\begin{equation*}
u(x,\tau)=\int_{-\infty}^{0}G(x,\tau,\xi,0)u^0(\xi)\ d\xi - \int_0^\tau M(\eta)G(x,\tau,s(\eta),\eta) \ d\eta+\int_0^\tau M(\eta)G(x,\tau,s_1(\eta),\eta)\ d\eta,
\end{equation*}
where $s_1(\tau)=s(\tau)+\dfrac{V_R}{\alpha(\tau)}$.
\qed
\begin{theorem}[Local existence of the non-linear Stefan-problem]
\label{SolMax_u}
Let $u^0(x)$ be a non-negative  function in
$C^0((-\infty,0])\cap C^1\left((-\infty, V_R)\cup(V_R,0]\right) 
\cap L^1((-\infty,0))$, such that $u^0(0)=0$ and suppose that 
$\frac{d u^0}{dx}$ vanishes at $-\infty$ and admits finite left and right 
limits at $V_R$. Then, there exists a unique maximal classical solution $u$ for
the problem \eqref{EqnEq} (equiva\-len\-tly \eqref{eq_u_definitiva}) in the sense of definition \ref{def_u}.
\end{theorem}
{\bf Proof.}
The system \eqref{EqnEq} (equivalently \eqref{eq_u_definitiva})
can be written as System  \eqref{AbstractEq} considering the function
$I(\omega)=\mu\left(-\log(\alpha(\omega))-D\right)\alpha(\omega)$,
on every time interval of size less 
than $\frac{1}{2}(e^{2D}-1)=\frac{\bar D}{2(1-\bar D)}$.
Then, for time $\tau \in [0,\frac{1}{2}(e^{2D}-1)]$ (equivalently
$\tau \in [0,\frac{\bar D}{2(1-\bar D)}]$), there is a unique
local solution defined on $[0,T_1]$ for System \eqref{EqnEq} (equivalently 
\eqref{eq_u_definitiva}), which is the solution for $I$ 
associated with the initial datum $M^0$.\\
If $T_1=\frac{1}{2}(e^{2D}-1)$, then we use $u(\cdot,T_1)$ and the values of $M$ in the interval $[0,\frac{1}{2}(e^{2D}-1)]$ as initial values for System \eqref{AbstractEq} in order to have a solution on $[0,T_2]$, $T_2\leq e^{2D}-1$.\\
We repeat this procedure until we find the maximal time of existence for 
the solution of \eqref{EqnEq} (equi\-va\-lent\-ly \eqref{eq_u_definitiva}).
\qed
The existence and uniqueness proved in 
Theorem \ref{SolMax_u} can be  translated into 
our initial system \eqref{delay_original}-\eqref{condiciones1}-\eqref{condiciones2} recovering  $\rho$ and $N$ by undoing the changes of variables 
\eqref{cambio1} and \eqref{cambio2}.

\begin{corollary}
	 Let $\rho^0$ be a non-negative $C^0((-\infty,V_F])\cap C^1((-\infty,V_R)\cup(V_R,V_F])\cap L^1((-\infty,V_F))$ function such that $\rho^0(V_F)=0$ and  $\dfrac{d \rho^0}{dv}$ decays at $-\infty$ and admits finite left and right limits at $V_R$. Then there exists a maximal $T^*\in(0,+\infty]$ and there exists a unique classical solution to the problem \eqref{delay_original}-\eqref{condiciones1}-\eqref{condiciones2} with $D>0$ on the time interval $[0,T^*)$.
\end{corollary}

\begin{remark}
	Using Duhamel's formula for $u$ and going back to $\rho$, we can see that as long as $\rho^0$ is fast decaying at $-\infty$ (for any polynomial function $f$, the quantities $f(v)\rho(v)$ and $f(v)\frac{d}{dv}\rho^0(v)$ go to 0 as $v$ goes to $-\infty$), then $\rho(\cdot,t)$ is fast decaying at $-\infty$ too for every positive $t$. This property will be implicitly used in other sections.
\end{remark}

\begin{remark}
	The same proof gives maximal classical solutions for the coupled excitatory-inhibitory system with positive delay studied in \cite{CS17-2}:
for a short time, the two equations decouple.
\end{remark}


\section{Global existence of solutions for the delayed model}

In this section we derive the main result of the work: the global existence 
of solutions for the delayed model \eqref{delay_original}.
The result is obtained directly for the average-inhibitory case
(as in the case without transmission delay \cite{CGGS}), while for the
average-excitatory case, it has to be derived through some of the properties 
of super-solutions.

\subsection{A criterion for the maximal time of existence}

The key step to obtain the main results of the paper is a criterion for 
the maximal time of existence of solutions, summarized in  Theorem 
\ref{th_maximal}. It ensures that solutions exists while the firing rate 
is finite.
With that purpose, first we show an auxiliary proposition, Proposition 
\ref{prop1}, which provides the tool to prove Theorem \ref{th_maximal}. 
Then using Theorem \ref{th_maximal} we derive Proposition \ref{prop2} 
which will allows to obtain the global existence of solutions for the 
inhibitory case. Their proofs are all omitted or sketched, since they are  
the same as in \cite{CGGS}[Proposition 4.1, Theorem 4.2, Proposition 4.3], 
because they are all consequences of the local existence result of Theorem 
\ref{SolMax_u}.

\begin{proposition}
Suppose that the hypotheses of Theorem \ref{SolMax_u} hold and that $u$ 
is a solution to \eqref{EqnEq}  (equivalently \eqref{eq_u_definitiva})
 in the time interval $[0,T]$. Assume in addition, that
\[
U_0:= \sup_{x \in (-\infty, s(t_0-\varepsilon)]} \abs{\partial_x u(x, t_0-\varepsilon)} < \infty \qquad \textrm{and that} \qquad  M^*= \sup_{t \in (t_0-\varepsilon, t_0)} M(t) < \infty,
\]
for some $0<\varepsilon<t_0\leq T$. Then,
$
\sup \{ \abs{\partial_xu(x,t)}: \ 
\ x \in (-\infty, s(t)], \ t \in [t_0-\varepsilon,t_0) \} < \infty,
$
with a bound depending only on the quantities $M^*$ and $U_0$.
\label{prop1}
\end{proposition}
Using this proposition we obtain the same criterion as in the case without
synaptic delay  \cite{CGGS}:

\begin{theorem}\label{th_maximal}
Suppose that the hypotheses of Theorem \ref{SolMax_u} hold. Then the solution $u$ can be extended up to a maximal time $0<\bar T \leq \infty$ given by
\[
\bar T= \sup \{ t\in (0,+\infty]: M(t)< \infty \}.
\]
\end{theorem}
In terms of the original system \eqref{delay_original}-\eqref{condiciones1}-\eqref{condiciones2}, we have the following  maximal time of existence result:

\begin{theorem}[Maximal time of existence]\label{resumen}
Let $\rho^0$ be a non-negative function in 
$C^0((-\infty,V_F])\cap C^1((-\infty,V_R)
\newline
\noindent
\cup(V_R,V_F])
\cap L^1((-\infty,V_F))$ 
such that $\rho^0(V_F)=0$ and  $\frac{d \rho^0}{dv}$ decays at $-\infty$ and admits finite left and right limits at $V_R$. Then there exists a unique maximal classical solution to the problem \eqref{delay_original}-\eqref{condiciones1}-\eqref{condiciones2} with $D\ge0$ on  the time interval $[0,T^*)$ where $T^*>0$ can be characterized by
\begin{equation*}
T^* =\sup \{t>0:N(t)<\infty \}.
\end{equation*}
\end{theorem}
{\bf Proof.}
The case $D=0$ is proved in the article \cite{CGGS}. In the case $D>0$, we use theorem \ref{th_maximal} and get the result directly.
\qed
Using Theorem \ref{th_maximal} we derive the key result for the global existence in the inhibitory case:
\begin{proposition}
Suppose that the hypotheses of Theorem \ref{SolMax_u} hold and that $u$ is 
a solution to  \eqref{eq_u_definitiva} (equivalently \eqref{EqnEq})
in the time interval $[0, t_0)$ for $b <0$. Then there exists 
$\varepsilon>0$ small enough and independent of $t_0$, such that, if
\begin{eqnarray}
 & \ & \bar U:=\sup_{x \in (-\infty, s(t_0-\varepsilon)]} \abs{\partial_x u(x, t_0-\varepsilon)} < \infty \label{cota44}
\end{eqnarray}
then $\sup_{t_0-\varepsilon<t<t_0} M(t) < \infty$, for $0<\varepsilon<t_0$.
\label{prop2}
\end{proposition}
Finally, combining Theorem \ref{th_maximal} with the previous result we 
obtain the global existence and uniqueness of classical solutions 
for the inhibitory case with synaptic delay for equivalent systems
 \eqref{EqnEq} and \eqref{eq_u_definitiva}.

\begin{proposition}\label{global_I}\label{global}
Suppose that the hypotheses of Theorem \ref{SolMax_u} hold and that $b<0$. 
Then there exists a unique global-in-time classical solution $u$ for system 
\eqref{eq_u_definitiva} (equivalently \eqref{EqnEq}) in the sense of 
Definition \ref{def_u} with initial datum $u^0$. 
Besides, if both $b$ and $b_0$ are negative, $s(t)$ is a 
monotone increasing function.
\end{proposition}
This proposition, translated to the initial delayed Fokker-Planck equation \eqref{delay_original} provides the global existence for the inhibitory case, 
as follows:
\begin{theorem}[Global existence - inhibitory case]\label{resumen2}
Let $\rho^0$ be a non-negative
  function 
in 
$C^0((-\infty,V_F])
\cap C^1((-\infty,V_R)\cup(V_R,V_F])\cap L^1((-\infty,V_F))$, such that $\rho^0(V_F)=0$ and $\partial_v\rho^0$ admits finite left and right limits at $V_R$. Suppose that  $\partial_v\rho^0$ decay at $-\infty$, then there exists a unique classical solution to the problem \eqref{delay_original}-\eqref{condiciones1}-\eqref{condiciones2} with $b<0$ and $D\ge0$ on the time interval $[0,T^*)$ with $T^* =+\infty$.
\end{theorem}

\subsection{Super-solutions and control over the firing rate}

We are not able  to obtain the global existence of solutions for the 
average-excitatory case  as it is done before for the average-inhibitory.
The difficulty is the extension
of the proposition \ref{prop2}
for the case $b>0$, which implies a uniform bound for $M$ in the 
average-excitatory
case.
Thus we have to proceed with  a different strategy, by means of a super-solution, to prove that the firing rate of any local solution cannot diverge in finite time. Then, applying the criterion of Theorem \ref{resumen}  the result is reached. We start introducing the notion of super-solution.

\begin{definition}
	Let  $T\in\R_+$, $D\ge0$ and $b_0=0$,
 ($\bar{\rho}$,$\bar{N}$) is said to be a (classical) super-solution to 
\eqref{delay_original}-\eqref{condiciones1}-\eqref{condiciones2} 
on $(-\infty,V_F]\times [0,T] $ if for all $t\in[0,T]$ we have 
$\bar{\rho}(V_F,t)=0$ and
	\begin{equation} \label{eq:upper}
\partial_t \bar{{\rho}} + \partial_v[(-v+b\bar{N}(t-D))\bar{{\rho}}]- a \partial_{vv}\bar{{\rho}}\geq  \delta_{v=V_R} \bar{N}(t), \quad \bar{N}(t)=- a \partial_v \bar{{\rho}}(V_F,t),
\end{equation}
	on $(-\infty,V_F]\times[0,T]$ in the distributional sense and on $((-\infty,V_F]\setminus V_R)\times [0,T]$ in the classical sense, with arbitrary values for $\bar{N}$ on $[-D,0)$.	
\end{definition}
We choose $b_0=0$ just for convenience as we can do it without loss of generality (as said in the introduction). Notice that for a solution in $C^{2,1}((-\infty,V_R)\cup(V_R,V_F]\times[0,T])\cap C^0((-\infty,V_R]\times [0,T])$, the condition reduces to satisfy the property in the classical sense in $(-\infty,V_R)\cup(V_R,V_F]\times[0,T]$ and having a decreasing jump discontinuity for the derivative on $V_R$ of size at least $\bar{N}/a$.

Notice also that when $T<D$, $\bar{N}(t-D)$ is an arbitrary initial datum, and thus if we find such a super-solution $\bar{\rho}$, then for every constant $\alpha>0$, the function $\alpha\bar{\rho}$ is also a super-solution.

We start proving the following comparison property between classical
solutions and super-solutions of \eqref{delay_original}-\eqref{condiciones1}-\eqref{condiciones2}. 

\begin{theorem}\label{thm:uppersol}
	Let $D>0$, $0<T<D$ and  $b_0=0$. 
Let $(\rho,N)$ be a classical solution of 
\eqref{delay_original}-\eqref{condiciones1}-\eqref{condiciones2} 
 on $(-\infty,V_F]\times[0,T]$ for the initial condition 
$(\rho^0,N^0)$ and let $(\bar{\rho},\bar{N})$ be a classical super-solution 
of \eqref{delay_original}-\eqref{condiciones1}-\eqref{condiciones2} 
on $(-\infty,V_F]\times[0,T]$. Assume that
	$$ \forall \, v\in(-\infty,V_F], \ \ \bar{\rho}(v,0) 
\ge \rho^0(v) \qquad \textrm{and} \qquad \forall \, t\in[-D,0),
\ \ \bar{N}(t)=N^0(t). $$
	Then,
	$$  \forall \, (v,t)\in (-\infty,V_F]\times[0,T], \ \ 
\bar{\rho}(v,t) \ge \rho(v,t) \qquad \textrm{and} \qquad \forall \, t\in[0,T], \ \  \bar{N}(t) \ge N(t). $$
\end{theorem}
{\bf Proof.} 
First, we prove that if $\bar{\rho}(v,t) \ge \rho(v,t)$ then 
$\bar{N}(t)\geqslant N(t)$. Due to the Dirichlet boundary 
condition for $\rho$ and the definition of super-solution we have 
$\rho(V_F,t)=\bar{\rho}(V_F,t)=0$ on $[0,T]$. 
Thus, as long as $\bar{\rho}(v,t) \ge \rho(v,t)$ holds, we have
$$ 
 -a\dfrac{\bar{\rho}(V_F,t)-\bar{\rho}(v,t)}{V_F - v} \ge -a\dfrac{\rho(V_F,t)-\rho(v,t)}{V_F - v}.  
$$ 
And taking the limit $v\to V_F$ we get
$\bar{N}(t)\geqslant N(t)$.

Then, denoting $w=\bar{\rho}-\rho$, we have for all $ (v,t)\in (-\infty,V_F]\times[0,T], $
$$ \partial_t w + \partial_v(-vw) + b\bar{N}(t-D)\partial_v\bar{\rho} - bN(t-D)\partial_v\rho - a \partial_{vv}w \ge  \delta_{v=V_R} (\bar{N}(t)-N(t)).   $$
As we assume $T<D$ we have by hypothesis $\bar{N}(t-D)=N^0(t-D)$ for all $t\in[0,T]$. Thus, as long as $w\ge 0$ holds, since $\bar{N}(t)\geqslant N(t)$,
$$ \partial_t w + \partial_v[(-v+bN^0(t))w] - a \partial_{vv}w \ge 0.$$
As $w(\cdot,0)\ge 0$, by a standard maximum principle theorem, we have
$\forall t\in[0,T], \ \ w(\cdot,t) \ge 0$,
and we conclude the proof.
\qed
Now, for fixed continuous (thus bounded) $N^0$ and the choice 
$\bar{N}(t)=N^0(t)$ in $[-D,0]$, we 
look for a super-solution  on $[0,D]$ of the form
\begin{equation}\label{rhobar}
\bar{\rho}(v,t)= e^{ \xi t} f(v),
\end{equation} 
where $\xi$ is large enough and $f$ is a carefully selected function, 
such that, the function $\bar{\rho}$ satisfies \eqref{eq:upper}, which means
\begin{equation}\label{Upper}
	(\xi-1)f + (-v+bN^0(t))  f' -a  f'' \geq  \delta_{v=V_R} V(t), \quad V(t)= - a f'(V_F). 
\end{equation}
 We show that $f$ defined as follow 
	$$  
	\begin{array}{rccl}
	f:& (-\infty,V_F] &\to      & \R_+ \\
	&     v         & \mapsto & \left\{\begin{array}{ll}
	1 & \ \mathrm{on} \ \left(-\infty,V_R\right]\\
	e^{V_R-v} \psi(v)+  \frac{1}{ \delta}(1-\psi(v)) (1-e^{\delta (v-V_F)})  & \ \mathrm{on} \ \left(V_R,V_F\right]
	\end{array}\right.
	\end{array}
	$$
verifies \eqref{Upper}.
To complete the definition of $f$ we explain which are $\psi$ and $\delta$:
\begin{enumerate}
	\item For $\varepsilon >0$ small enough, such that $\frac{V_F+V_R}{2}
        + \varepsilon <V_F$, we consider  $\psi \in C^{\infty}_b(\R)$ satisfying
	$0 \leq \psi \leq 1$ and
	$$
        \psi \equiv 1 \hbox{ on } \left(-\infty, \frac{V_F+V_R}{2}\right) 
        \hbox{ and } \psi \equiv 0 \hbox{ on }  \left(\frac{V_F+V_R}{2} + 
        \varepsilon, +\infty\right).
        $$
	\item For $B>0$, such that
        $\left|-v+bN^0(t)\right| \leq B$, $\forall \, t \in [-D,0)$, 
        $\forall \, v \in \left(V_R,V_F\right)$,
	we take $\delta >0$ such that $a \delta - B \geq 0$.
\end{enumerate}
Notice that $f$ being the sum of two continuous non-negative functions that never vanish at the same point, we have
$$\inf_{v\in\left(V_R,\frac{V_F+V_R}{2}+ \varepsilon\right)}f(v) > 0.
$$
With these choices, $\bar{\rho}(v,t)$ is a super-solution on $[0,D]$ for 
$\xi$ large enough, because:
\begin{itemize}
	\item On $\left(-\infty,V_R\right)$, $\bar{\rho}$ is independent of 
        $v$, thus the definition is satisfied if and only if $\xi > 1$.
	\item Around the $V_R$ point the inequality \eqref{Upper} has to hold 
        in the sense of distribution, that is in our case
	$$f'(V_R^{+}) - f'(V_R^{-}) \leqslant f'(V_F) $$
	This inequality is satisfied because $f'(V_R^{-})=0$, 
        $f'(V_R^{+})=-1$ and $f'(V_F)=-1$.
	\item On $\left(V_R,\frac{V_F+V_R}{2}+ \varepsilon\right)$, 
        we choose $\xi$ such that 
	$$  
        (\xi-1) \inf_{v\in\left(V_R,\frac{V_F+V_R}{2}+ \varepsilon\right)}f(v) 
        \geq \sup_{v\in \left(V_R,\frac{V_F+V_R}{2}+ \varepsilon\right)}
        \left(B| f'(v)| +a |f''(v)|\right), 
        $$
	which is possible because $\inf_{v\in\left(V_R,\frac{V_F+V_R}{2}+ 
        \varepsilon\right)}f(v) > 0$. Then the super-solution inequality 
        \eqref{Upper} holds.
	\item  On $\left(\frac{V_F+V_R}{2}+ \varepsilon,V_F\right)$, 
        the desired inequality \eqref{Upper} holds since
	$$ (-v+b N^0(t)) f' -a f'' = 
           e^{\delta(V_F-v)}\left[a\delta-(-v+b N^0(t))\right] 
           \geq e^{\delta(V_F-v)}\left[a\delta-B\right] \geq 0.$$
\end{itemize}
Given this super-solution on $[0,D]$ for any fixed continuous $N^0(t)$, 
we can prove global existence for local solutions.

\begin{theorem}[Global existence - excitatory and inhibitory cases] \label{thm:globalexist}
	Let $\rho^0$ be a non-negative function
in $C^0((-\infty,V_F])\cap C^1((-\infty,V_R)\cup(V_R,V_F])\cap L^1((-\infty,V_F))$,  such that $\rho^0(V_F)=0$ and  $\frac{d \rho^0}{dv}$ decays at $-\infty$ and admits finite left and right limits at $V_R$; let $N^0\in C^0([-D,0])$. Let $(\rho,N)$ be the corresponding maximal classical solution of \eqref{delay_original}-\eqref{condiciones1}-\eqref{condiciones2} with $D>0$. Then, the maximal existence time for the local solution $(\rho,N)$ is $T^*=+\infty$.
\end{theorem}
{\bf Proof.}
Assume the maximal time of existence $T^*$ is finite, this means,
using Theorem \ref{resumen}, that the firing rate $N$ diverges 
when $t\to T^*$. We prove that this is a contradiction with the
fact that $\bar{\rho}$ given by \eqref{rhobar} is
a super-solution.

 As the maximal solution was showed previously to be unique, we assume without loss of generality that $T^*=\frac{D}{2}<D$ by using the new initial conditions
$$\tilde{\rho}^0(v)=\rho(v , T^*-\frac{D}{2}\ ) \quad \forall \ v \in(-\infty, V_F] \quad \mathrm{and} \quad \tilde{N}^0(\tilde t)=N\left(T^*-\frac{D}{2} +\tilde t\right), \quad \tilde t\in[-D,0).$$
As $\tilde \rho^0$  is continuous and vanish at $V_F$ and $-\infty$, it belongs to $L^\infty((-\infty,V_F])$ and therefore there exists $\alpha\in\R_+^*$ such that the super-solution $\bar{\rho}$ we constructed satisfies
$\alpha\bar{\rho}(v,0) \ge \tilde\rho^0(v)$, for all $v\in(-\infty,V_F]$,
where we use the fact that $\bar{\rho}$ never vanish on $(-\infty,V_F)$.
Then, by Theorem \ref{thm:uppersol}, we have
$$   N\left(T^*-\frac{D}{2} +\tilde t\right)=\tilde N(\tilde t)\leq \bar{N}(\tilde t) = a e^{\xi \tilde t} \quad \forall \, \tilde t\in \left[0,\frac{D}{2}\right). $$
Thus,
$N(t)\leq a e^{\xi(t-T^*+\frac{D}{2})}$ for all $t \in 
\left[T^*-\frac{D}{2}, T^*\right)$.
\noindent Therefore, by continuity, there is no divergence of the firing rate $N$ when $t\rightarrow T^*$, and thus by Theorem \ref{resumen} we reach a contradiction. 
\qed

\begin{remark}
	As we remark for the local existence, the result about global
 existence works for the coupled excitatory-inhibitory system with positive 
delay studied in \cite{CS17-2}. 
\end{remark}

\section{Qualitative properties and long time behavior of solutions}

The aim of this part is to extend the results about long time behavior obtained in  \cite{CCP} and \cite{carrillo2014qualitative}  for the case without delay to the case where a synaptic delay is considered. All the results can be extended, assuming that the delay is small enough with respect to  the other parameters of the model, in particular in comparison with the connectivity $b$.


In all this section, we use two standard technical results fully stated in
 Appendix \ref{TeXniK} (entropy  method and  Poincar\'{e}'s inequality). 
We also recall here that a steady state of the system 
\eqref{delay_original}-\eqref{condiciones1}-\eqref{condiciones2} is defined 
as a solution of the problem
\[
\dfrac{\partial}{\partial v}\left[(-v + bN_\infty)\rho_{\infty}\right]-a\dfrac{\partial^2 \rho_{\infty}}{\partial v^2}= \ \delta_{v=V_R}  N_\infty, \label{NLIF}	
\]
\[
N_\infty(t)= -a \dfrac{\partial \rho_{\infty}}{\partial v}(V_F),  \quad \rho_\infty(V_F)=0, \quad \rho_{\infty}(-\infty)=0,
\]
\[
\rho_\infty(v) \geqslant 0, \quad \int_{-\infty}^{V_F} \rho_{\infty}(v)d v = 1.
\]
Such steady states are exactly the same for the case without 
transmission delay ($D=0$) and for the case with synaptic delay ($D>0$). 
They are studied in \cite{CCP} and an implicit form is given:
\[  \rho_{\infty}(v) = \dfrac{N_\infty}{a}e^{-\frac{(v-bN_\infty)^2}{2a}}   \int_{\max(v,V_R)}^{V_F}  e^{\frac{(w-bN_\infty)^2}{2a}} dw.  \]
Moreover, it is proved in \cite{CCP} that for small enough positive $b$ or every negative $b$, there is an unique steady state $(\rho_{\infty},N_\infty)$. For $b>0$ large enough, there are no steady states.\\

Using the results of \cite{CCP}, it is also possible to directly deduce for small enough values of $b$ (\textit{via} the continuous dependence in $b$ and the value at 0) that
\begin{equation}\label{limb}
\lim_{b\to0} N_\infty(b) > 0,
\end{equation}
where we denote by $N_\infty(b)$ the stationary firing rate for the system
\eqref{delay_original}-\eqref{condiciones1}-\eqref{condiciones2}
with connectivity parameter $b$.
\subsection{Uniform estimates on the firing rate}

We  separate the average-inhibitory  and the average-excitatory cases. Indeed, in the average-inhibitory case $(b<0)$, all the results are uniform with respect to the initial data, which is not the case when $b>0$ (average-excitatory case), in accordance with the blow-up structure of Equation \eqref{delay_original} in the case without delay. In the previous section we proved the global existence of the solution when synaptic delay is considered. However, in the proofs, the structure of Equation \eqref{delay_original}, in the average-excitatory case, leads to results which strongly depends on the initial data.\\\\
The following theorem, which is similar to Theorem 3.1 of \cite{carrillo2014qualitative}, provides some $L^2$ control over the firing rates, and thus, supplies the  main tool to prove the long time behavior result of Theorem \ref{thm:expdecay}.

\begin{theorem}\label{L24}
	Let $b_1>0$ such that there exists a stationary state of System \eqref{delay_original}-\eqref{condiciones1}-\eqref{condiciones2} and let $\rho_\infty^1$ be the corresponding  stationary state. Let $V_M$ with  $V_F>V_M>V_R$ and define
	$
	S(b_1,V_M):= \int_{V_M}^{V_F} \frac{\left(\rho^0\right)^2}{ \rho_\infty^1}\ dv
	$,  
	with an initial datum chosen such that $S(b_1,V_M)<+\infty.$ Then:
	\begin{enumerate}
		\item[i)] There exists a constant $C$ independent of $S(b_1,V_M)$ and there exists  a time $T>0$ depending only on $V_M$ and $S(b_1,V_M)$, such that for all intervals $J \subset( T,+\infty)$ and for all $b \leq 0$:
		\begin{equation}\label{contrN1}
		\int_{J} N(t)^2\ dt \leq C(1+|J|).
		\end{equation}
		
		\item[ii)] Assume $b>0$ is small enough depending on $S(b_1,V_M)$, the delay $D$ and $V_M$, then there exists a constant $C$ such that for all $t \geq 0$:
		\begin{equation}\label{contrN2}
		\int_{0}^t N(t)^2\ dt \leq C(1+t).
		\end{equation}
		
	\end{enumerate}
\end{theorem}
\begin{remark}
	Let us mention that the condition on the delay only appear in the 
average-excitatory case in Theorem \ref{L24}.
\end{remark}
{\bf Proof.}
The main idea of the proof is to use entropy type inequalities and to compare the solution of Equation \eqref{delay_original} with a kind of super-solution which is a stationary state of Equation \eqref{delay_original} with a strictly positive connectivity parameter in a similar manner as in \cite{carrillo2014qualitative} [Theorem 3.1], with the difficulty that now there is a synaptic delay in the 
equation. 

Let us  first introduce some notations and a useful function. We set  the function $\gamma$ as in \cite{carrillo2014qualitative} and defined by
$$ \forall \, v\in(-\infty,V_F], \ \ \gamma (v)= 1_{ v > \alpha }e^{\frac{-1}{\beta-(V_F-v)^2}},$$
where $\beta=(V_F-\alpha)^2$ and $\alpha\in(-\infty,V_F)$. Given a steady state $(\rho_\infty^1,N_\infty^1)$ associated to the parameter $b_1$ and a solution $(\rho,N)$ associated to the parameter $b$ and the initial condition $(\rho^0,N^0)$, we denote in the following
$$ H(v,t)=\frac{\rho(v,t)}{\rho_\infty^1(v)} \qquad \mathrm{and} \qquad W(v,t)=\rho_\infty^1(v)(H(v,t))^2 = \frac{\rho(v,t)^2}{\rho_\infty^1(v)}$$ 
We also assume $V_R$ positive without loss of generality (remember \eqref{traslado}).

Let us first prove Theorem \ref{L24} in the inhibitory case $( b \leq 0)$. Let  
$I(t)=\int_{-\infty}^{V_F}W(v,t)\gamma(v)\ dv$.  Following the computations made in \cite{carrillo2014qualitative}[Section 3.1.1], we find  that, for all $\alpha$ close enough to $V_F$, there exist $C_1,C_2,C_3\in\R_+^*$ independent of $b$ such that
\begin{equation}\label{integramos}
 \frac{d}{dt}\int_{-\infty}^{V_F}W(v,t)\gamma(v)\ dv \leq -C_1N(t)^2+C_2 +(bN(t-D)-C_3)\int_{-\infty}^{V_F}W(v,t)\gamma(v)\ dv.
\end{equation} 
We have, due to $I(0)<+\infty$ and $b\leq 0$,
$$ \frac{d}{dt} I(t) \leq C_2-C_3I(t).$$
Thus, there exists a time $T>0$ depending on $I(0)$, such that,
for  $t>T$:
$I(t)\leq 2\frac{C_2}{C_3}$.
This conclude the proof of the first part of Theorem \ref{L24}   by integrating the differential inequality \eqref{integramos} on every interval $I\subset[T,+\infty[$.

Let us now assume that $b>0$. Let us again consider $I(t)=\int_{-\infty}^{V_F}W(v,t)\gamma(v)\ dv$. Following again the computations made in \cite{carrillo2014qualitative}[Section 3.1.2], we obtain that there exist $C_1$, $C_2$, $C_3$, $C_4$ and $C_5=\frac{\gamma(V_F)}{N_\infty^1}$ depending only on $\gamma$ and $\rho_\infty^1$ such that for all $t>0$,
\begin{equation}\label{Ayda}
\frac{d}{dt}I (t)\leq \Big(-C_4 + C_2b^2N(t-D)^2\Big)I(t) + C_1b^2 N(t-D)^2 - C_5 N(t)^2 + C_3.
\end{equation} 
Let $\bar{M}=\max({I(0),\frac{C_3}{C_4}})$, where we assume $I(0)<\infty$. As $I(0)$ is strictly bounded by $2\bar M$, by continuity of $I$, there exists a closed interval $[0,\eta]$ on which $I$ is bounded by $2\bar M$. Thus, let $[0,\bar B]$ ($[0,\bar B)$ if $\bar B=+\infty$)   the largest closed interval on which $I$ is bounded by $2\bar M$. Assume $\bar B$ is finite.
On this interval $[0,\bar B]$, we have 
$$ \frac{d}{dt} I(t)\leq -C_4I(t) + \Big(C_1+2C_2\bar M\Big)b^2 N(t-D)^2 - C_5 N(t)^2 + C_3. $$
Let $C_6=C_1+2C_2\bar M$. As $C_3 - C_4\bar M \leq 0$, we have
$$\frac{d}{dt}(I(t)-\bar M)\leq -C_4(I(t)-\bar M) + C_6b^2 N(t-D)^2 - C_5 N(t)^2. $$
Gronwall's Lemma implies that, for all $t\in[0,\eta[$,
\begin{equation}  I(t)-\bar M \leq e^{-C_4t}(I(0)-\bar M) +e^{-C_4 t}\int_{0}^{t} \Big(C_6b^2N(s-D)^2 - C_5N(s)^2 \Big)e^{C_4s}\ ds. \end{equation}
Therefore, we can write 
\begin{align*}
  \ \int_{0}^{t} \Big(C_6b^2N(s-D)^2 - C_5N(s)^2 \Big)e^{C_4s}\ ds 
 = & \ C_6b^2\int_{0}^{t} N(s-D)^2e^{C_4s}\ ds - C_5\int_{0}^{t}N(s)^2e^{C_4s}\ ds\\
 \leq &  \ C_6b^2\int_{-D}^{0} N(s)^2e^{C_4s + C_4D}\ ds\\
= & \ C_6e^{C_4D}b^2\int_{-D}^{0} N(s)^2e^{C_4s}\ ds\\
 \leq & \ C_6e^{C_4D}b^2\int_{-D}^{0} N^0(s)^2\ ds.
\end{align*}
Hence,
$I(t) \leq \bar M + e^{-C_4t}C_6e^{C_4D}b^2\int_{-D}^{0} N^0(s)^2\ ds$.
For $b$ small enough,
$$ \bar M + e^{-C_4t}C_6 e^{C_4D}b^2\int_{-D}^{0} N^0(s)^2\ ds  < 2 \bar M,$$
which implies $I(\bar B)<2\bar M$ and by continuity of $I$ there exists $\varepsilon\in\R_+^*$ such that $I$ is stricly bounded by $2 \bar M$ on $[\bar B,\bar B+\varepsilon]$. This is a contradiction to the maximality of $\bar B$. Thus, $\bar B=+\infty$ and $I(t)\leq 2 \bar M$, for all
$t\in\R_+$.
Hence, integrating \eqref{Ayda} between 0 and $t$ we obtain
$$ -I(0) \leq  \Big(C_1 + 2C_2\bar M \Big)b^2\int_{-D}^{t-D}N(s)^2\ ds -C_5\int_{0}^{t}N(s)^2\ ds + C_3 t,$$
that gives, for $b$ small enough, 
$$  \frac{C_5}{2}\int_{0}^{t}N(s)^2\ ds \leq I(0) + \Big(C_1 + 2C_2\bar M \Big)b^2\int_{-D}^{0}N^0(s)^2\ ds + C_3t,$$
which finishes the proof of Theorem \ref{L24}.
\qed 
\subsection{Convergence toward steady states for small connectivities}
Now, we can prove our main result of this section, that ensures that even with delay, the solutions converge exponentially fast to the steady state, when the connectivity parameter is small, and with an appropriate initial condition.

\begin{theorem}\label{thm:expdecay}
	Let $b_1$ such that there exists a steady state $(\rho_\infty^1,N_\infty^1)$ for \eqref{delay_original}-\eqref{condiciones1}-\eqref{condiciones2} with $b_0=0$. Let $V_M\in(V_R,V_F)$, $(\rho^0,N^0)$ an initial condition and
	$S(b_1,V_M):=\int_{V_M}^{V_F}\frac{\rho^0(v)^2}{\rho_\infty^1(v)}\ dv$.
	If $S(b_1,V_M)<+\infty$, and if there exists $C_0$ such that $\rho^0(\cdot) \leq C_0 \rho_{\infty}^1(\cdot)$, then: 
	\begin{itemize}
		\item  If $b>0$, there exists $\eta\in\R_+^*$ depending only on $\rho^0$, $N^0$, $S(b_1,V_M)$, $D$ and $V_M$ such that if $0<b<\eta$, then there exist $A_0,\mu,\gamma \in\R_+^*$ depending only on $(\rho^0,N^0)$ such that
\begin{align*}
\int_{-\infty}^{V_F}\rho_\infty\left( \frac{\rho-\rho_\infty}{\rho_\infty} \right)^2(v,t)\ dv \leq & A_0e^{-\mu t}\Bigg[\int_{-\infty}^{V_F}\rho_\infty\left( \frac{\rho^0-\rho_\infty}{\rho_\infty} \right)^2(v)\ dv  \\
& +\ \frac{8b^2}{a} \int_{-D}^{0}\left( N^0(s) - N_\infty \right)^2e^{\gamma (s+D) - \int_{-D}^{0}(N^0(u)-N_\infty)^2\ du }\ ds \Bigg].
\end{align*}
		\item  If $b< 0$, there exists $\eta\in\R_+^*$, independent of $S(b_1,V_M)$ such that if $-\eta < b < 0$, there exist $T,\mu,\gamma \in\R_+^*$ depending only on the initial condition, and $B_0\in\R_+^*$ depending only on the values of $N$ on $[0,T]$, such that for all $t\in(T,+\infty)$,
\begin{align*}
 \int_{-\infty}^{V_F}\rho_\infty\left( \frac{\rho-\rho_\infty}{\rho_\infty} \right)^2(v,t)\ dv \leq &
		B_0e^{-\mu t}\Bigg[\int_{-\infty}^{V_F}\rho_\infty\left( \frac{\rho^0-\rho_\infty}{\rho_\infty} \right)^2(v)\ dv \\
& +\ \frac{8b^2}{a} \int_{-D}^{0}\left( N^0(s) - N_\infty \right)^2e^{\gamma (s+D) - \int_{-D}^{0}(N(u)-N_\infty)^2\ du }\ ds \Bigg]  .
\end{align*}
	\end{itemize}
\end{theorem}
{\bf Proof.}
Let  $G:\R_+\to\R$  be a convex function of class $C^2$. As $\rho^0\leq C_0 \rho_{\infty}^1$ and because the rate of decay for $\rho_{\infty}(b)$ is asymptotically the same for every $b$ at $-\infty$, we have for every $b$ such that a stationary state exists another constant $\tilde{C_0}$ such that $\rho^0(\cdot)\leq\tilde{C_0}\rho_{\infty}(b)(\cdot)$. Then, the following entropy equality holds, applying theorem \ref{thm:entropy}.
\begin{align*}
 \  \frac{d}{dt}\int_{-\infty}^{V_F} \rho_{\infty}
G\left(\frac{\rho}{\rho_{\infty}}\right) \ dv   
= & \  -a\int_{-\infty}^{V_F}\rho_\infty \left[\frac{\partial}{\partial v}\left(\frac{\rho}{\rho_\infty}\right)\right]^2G''(\frac{\rho}{\rho_{\infty}})\ dv   \\
& \  -N_{\infty}\left[G\left(\frac{N}{N_\infty}\right) - G\left(\frac{\rho}{\rho_\infty}\right) - \Big( \frac{N}{N_\infty} - \frac{\rho}{\rho_\infty} \Big)G'\left(\frac{\rho}{\rho_\infty}\right)   \right]\bigg\rvert_{v=V_R}  \\
& \   + b(N(t-D)-N_\infty)\int_{-\infty}^{V_F}\frac{\partial \rho_\infty}{\partial v}\left[G\left(\frac{\rho}{\rho_\infty}\right) -\frac{\rho}{\rho_\infty}G'\left(\frac{\rho}{\rho_{\infty}}\right)\right]\ dv.
\end{align*}
Now, denote  
$q(v,t)=\frac{\rho(v,t)}{\rho_\infty(v)}$ and $F(t)=\frac{N(t)}{N_\infty}$.
We choose $G(x)=(x-1)^2$. Hence, for all $\varepsilon\in(0,\frac{1}{2})$,  using the inequality $(a+b)^2\ge \varepsilon(a^2-2b^2)$),
 yields
$$  G(F(t))-G(q(V_R,t))-\Big(F(t)-q(V_R,t)\Big)G'(q(V_R,t)) = (F(t)-q(V_R,t))^2
\geq \varepsilon(M(t)-1)^2 - 2\varepsilon(q(V_R,t)-1)^2.  $$
Then, according to the expression of $\rho_\infty$, the Poincar\'{e}-like inequality (see Appendix \ref{TeXniK}) and the Sobolev injection of $L^\infty(I)$ in $H^1(I)$ for a sufficiently small neighborhood $I$ of $V_R$, there exists $C\in\R_+^*$ such that
$$ \left| q(V_R,t) - 1  \right|^2 \leq C\int_{-\infty}^{V_F} \rho_\infty(v)\left(\frac{\partial q}{\partial v}\right)^2(v,t)\ dv.$$
Thus, for $\varepsilon$ satisfying $2CN_\infty\varepsilon\leq \frac{a}{2}$, we have
\begin{eqnarray}\label{Maj1}
 & \ & -N_{\infty}\left[G\left(\frac{N}{N_\infty}\right) - G\left(\frac{\rho}{\rho_\infty}\right) - \Big( \frac{N}{N_\infty} - \frac{\rho}{\rho_\infty} \Big)
G'\left(\frac{\rho}{\rho_{\infty}}\right)   \right]\bigg\rvert_{v=V_R} \nonumber \\
& \leq & -N_\infty \varepsilon G\left(\frac{N(t)}{N_\infty}\right)+\frac{a}{2}\int_{-\infty}^{V_F}\rho_\infty\left(\frac{\partial q}{\partial v}\right)^2(v,t)\ dv.
\end{eqnarray}
On the other hand, as $G(x)-xG'(x)=1-x^2$, we have by integration by parts,
\begin{align*}
 b(N(t-D)-N_\infty)\int_{-\infty}^{V_F}\frac{\partial \rho_\infty}{\partial v}\Big[G(q) -qG'(q)\Big]\ dv 
=   2b(N(t-D)-N_\infty)\int_{-\infty}^{V_F} \rho_\infty q\frac{\partial q}{\partial v}\ dv.
\end{align*}
Using the inequality $cd\leq \bar\varepsilon c^2 + \frac{1}{\bar\varepsilon}d^2$ with $\bar\varepsilon=\frac{a}{2}$, $c=\sqrt{\rho_\infty}\frac{\partial q}{\partial v}$ and $d=2b(N(t-D)-N_\infty)\sqrt{\rho_\infty}q$, we obtain 
\begin{eqnarray}\label{Maj2}
 & \ & b(N(t-D)-N_\infty)\int_{-\infty}^{V_F}\frac{\partial \rho_\infty}{\partial v}\Big[G(q) -qG'(q)\Big]\ dv \nonumber \\
& \leq & \  \frac{a}{2}\int_{-\infty}^{V_F} \rho_\infty \left(\frac{\partial q}{\partial v}\right)^2\ dv + \frac{8(b(N(t-D)-N_\infty))^2}{a}\int_{-\infty}^{V_F} \rho_\infty q^2\ dv \nonumber\\
& \leq & \frac{a}{2}\int_{-\infty}^{V_F} \rho_\infty \left(\frac{\partial q}{\partial v}\right)^2\ dv + \frac{8(b(N(t-D)-N_\infty))^2}{a}\int_{-\infty}^{V_F} \rho_\infty (q-1)^2\ dv\\
& \ & +\frac{8(b(N(t-D)-N_\infty))^2}{a}.\nonumber
\end{eqnarray}
Collecting the previous bounds \eqref{Maj1} and \eqref{Maj2}, we have
\begin{align*}
 \frac{d}{dt}\int_{-\infty}^{V_F}\rho_\infty(v)G(q(v,t))\ dv  
\leq & \  -N_\infty \varepsilon G\left(\frac{N(t)}{N_\infty}\right) + 
\frac{8(b(N(t-D)-N_\infty))^2}{a}\int_{-\infty}^{V_F}\rho_\infty(v)
G(q(v,t))\ dv \\
& \  -a\int_{-\infty}^{V_F} \rho_\infty\left(\frac{\partial q}{\partial v}\right)^2\ dv +\frac{8(b(N(t-D)-N_\infty))^2}{a}.
\end{align*}
And using the Poincar\'{e}'s inequality 
and the specific form of $G$ we chose, we obtain 
\begin{align*}
  \frac{d}{dt} \int_{-\infty}^{V_F}\rho_\infty(v)G(q(v,t))\ dv 
\leq & \  -\frac{\varepsilon}{N_\infty}\Big(N(t)-N_\infty\Big)^2+\frac{8b^2}{a}\Big(N(t-D)-N_\infty\Big)^2  \\
& \  +\left( -\nu a + \frac{8(b(N(t-D)-N_\infty))^2}{a}\right)\int_{-\infty}^{V_F}\rho_\infty(v)G(q(v,t))\ dv.  
\end{align*}
Then, denoting
$E(t)=\int_{-\infty}^{V_F}\rho_\infty(v)G(q(v,t))\ dv$,
and
$\phi(t)= -\nu at + \frac{8b^2}{a}\int_{0}^{t}(N(s-D)-N_\infty)^2\ ds$,
and using Gronwall's lemma, the previous inequality becomes 
\begin{align*}
  E(t)\leq & \ e^{\phi(t)}E(0) + \frac{8b^2}{a}e^{\phi(t)}\int_{0}^{t}(N(s-D)-N_\infty)^2e^{-\phi(s)}\ ds \\
& \ -\frac{\varepsilon}{N_\infty}e^{\phi(t)}\int_{0}^{t}(N(s)-N_\infty)^2e^{-\phi(s)}\ ds\\
\leq & \  e^{\phi(t)}E(0) + \frac{8b^2}{a}e^{\phi(t)}\int_{-D}^{t-D}(N(s)-N_\infty)^2e^{-\phi(s+D)}\ ds\\
& \ -\frac{\varepsilon}{N_\infty}e^{\phi(t)}\int_{0}^{t}(N(s)-N_\infty)^2e^{-\phi(s)}\ ds.
\end{align*}
But, for $s>0$,
$$-\phi(s+D) = \nu as +\nu aD - \frac{8b^2}{a}\int_{0}^{s+D}(N(u-D)-N_\infty)^2\ du$$
$$	\leq -\phi(s) +\nu aD.$$
Therefore, as using \eqref{limb} we can choose $b$ such that
$\frac{8b^2 e^{\nu a D}}{a}\leq \frac{\varepsilon}{N_\infty(b)}$,
we can write
\begin{align*}
  E(t) \leq & \  e^{\phi(t)}E(0) + \frac{8b^2}{a}
e^{\phi(t)}
\int_{-D}^{0}(N(s)-N_\infty)^2e^{-\phi(s+D)}\ ds \\
& \ -\frac{\varepsilon}{N_\infty}e^{\phi(t)}\int_{t-D}^{t}(N(s)-N_\infty)^2e^{-\phi(s)}\ ds \\
\leq& \  e^{\phi(t)}\Bigg( E(0) + \frac{8b^2}{a}\int_{-D}^{0}(N^0(s)-N_\infty)^2e^{-\phi(s+D)}\ ds\Bigg). 
\end{align*}
\begin{itemize}
\item If $b>0$ :
By Theorem \ref{L24}, in the average-excitatory case, there exist constants $\beta_1,\beta_2\in\R_+$ such that
$$ \phi(t) \leq -\nu at + \frac{8b^2}{a}(\beta_1 + \beta_2 t).$$
Thus, for $b$ small enough ($\beta_1(b)$ and $\beta_2(b)$ do not grow when $b$ goes to 0), $e^{\phi(t)} \leq A_0e^{-\mu t}$,
with $A_0$, $\mu$ positive real numbers.
\item If $b\leq0$: 
By Theorem \ref{L24}, in the inhibitory case, there exists $T\in\R_+$ and there exist constants $\beta_1,\beta_2\in\R_+$ such that
$\phi(t) \leq -\nu at + \frac{8b^2}{a}(\beta_1 + \beta_2 (t-T))$, for 
$t\geq T$.
Thus, for $b$ small enough (these $\beta_1$ and $\beta_2$ don't depend on $b$) and $t$ large enough,
$e^{\phi(t)} \leq B_0e^{-\mu t}$,
with $B_0$, $\mu$ positive real numbers and $B_0$ depending of the values of $N$ between $0$ and $T$. \qed
\end{itemize}

\subsection{Non-existence of periodic solutions for a large connectivity parameter}

Numerical results for the NNLIF models  have been presented in \cite{CS17-2}. 
On one hand, we observe how the blow-up is avoided if we include a synaptic delay. For this case, we have a small value of $b$ combined with a concentrated initial condition, which produces the blow-up of the solution without delay \cite{CCP}. For this value of $b$ there is a unique steady state \cite{CCP}, and the solution seems to tend to it, after avoiding the blow-up due to the delay.

On the other hand, also it is known that a blow-up situation happens  for 
large value of $b$ \cite{CCP}. If we include the delay, the solutions avoid the blow-up, but they do not tend to an equilibrium, since for large values of $b$ there is no steady state \cite{CCP}. Numerically, 
the firing rates  seem to grow slowly all the time with limit $+\infty$, but without diverging in finite time \cite{CS17-2}. 
In this case it could be  expected solutions to present a somehow 
periodic behavior, but  it did not observed numerically.

Here we clarify a bit the situation by proving analytically that it is impossible for periodic solutions to exist when $b$ is over $V_F-V_R$ in the case $V_F\leq 0$ and $b_0=0$.

\begin{theorem}
	If $b>V_F-V_R$, $b_0=0$ and $V_F\leq 0$, then for any $D\geq 0$ there are no classical periodic solutions to equation \eqref{delay_original}-\eqref{condiciones1}-\eqref{condiciones2} such that
	$\int_{-\infty}^{V_F}|v|\rho^0(v)d v < +\infty$.
\end{theorem}
{\bf Proof.}
Assume there exists a $T-$periodic solution $(\rho,N)$ of equation
$$\partial_t \rho + \partial_v\big[ (-v+bN(t-D))\rho\big]-a \partial_{vv}\rho=  \delta_{V_R}(v) N(t)$$
such that $v\rho^0\in L^1((-\infty,V_F))$. Then, as we said before, $v\mapsto v\rho(v,t)\in L^1((-\infty,V_F))$ for all time since the decay assumption propagates.
Denoting
$\Phi(v):= \frac{1}{T} \int_{0}^{T} \rho(v,t) dt$,
we have
$\int_{-\infty}^{V_F} \Phi(v)dv=1$
and,   $\Phi$ satisfies
\begin{equation}\label{eq1}
\partial_v\big( -v \Phi + b  \frac{1}{T}\int_{0}^{T} N(t-D) \rho(v,t)dt\big) -a \partial_{vv} \Phi = \delta_{v=V_R} \overline{N},
\end{equation}
where   $\overline{N}:= \frac{1}{T} \int_{0}^{T} N(t)dt=  \frac{1}{T} \int_{0}^{T} N(t-D)dt$.
Now, we multiply equation \eqref{eq1} by $v$ and we integrate: 
$$   \int_{-\infty}^{V_F}  v \Phi dv -b \overline{N} + (V_F-V_R)  \overline{N}=0.    $$
Hence,
$  \int_{-\infty}^{V_F}  v \Phi dv=  \big(b- (V_F-V_R)\big) \overline{N}$. 
As $V_F \leq 0$, the integral is negative, leading to $(b- (V_F-V_R)) <0$,
which is a contradiction.
\qed

\begin{corollary}
	If $b>V_F-V_R$, $V_F\leq 0$ and $b_0=0$, there is no steady state to equation \eqref{delay_original}-\eqref{condiciones1}-\eqref{condiciones2}.
\end{corollary}

\begin{remark}
	This corollary upgrade a result of \cite{CCP} that classifies the number of steady states. Indeed, the theorem 3.1 of \cite{CCP} tell us (in the case $b_0=0$) that there is no steady state when
		\[  b > \max\left( \  2(V_F-V_R)\ ,\ 2V_F \int_{0}^{+\infty} e^{-\frac{s^2}{2}} \dfrac{e^{\frac{s V_F}{\sqrt{a_0}} } - e^{\frac{s V_R}{\sqrt{a_0}} } }{s}d s\  \right), \]
	and that there are at least two steady states when
$b > V_F - V_R$  and $0 < 2a_0 b < (V_F-V_R)^2V_R$. 
	The case $b>V_F-V_R$ and $V_F\leq 0$ was not covered by the result.
\end{remark}

\section{Conclusions and open problems}

This paper focuses on the NNLIF system with a synaptic delay.  For the
average-excitatory case without transmission delay, solutions can
blow-up in finite time \cite{CCP} due to the divergence of the firing
rate \cite{CGGS}.  Nevertheless, at microscopic level, it has been
proved that if the synaptic delay is taken into account, solutions are
always global-in-time \cite{delarue2015particle}. Moreover, for
delayed NNLIF system numerical results show that the blow-up
phenomenon is avoided \cite{CS17-2}. In this work we prove the
global-in-time existence for the delayed NNLIF model, confirming the
previous numerical observations. The techniques developed in
\cite{CGGS}, for the case without synaptic delay, are not enough for
the average-excitatory delayed system, since it is not possible to
find a uniform bound for the firing rate. We reach the proof by combining
these techniques with the construction of super-solutions
\cite{carrillo2014qualitative}, which provide a control of the firing
rate.  Moreover, we show qualitative properties of the solutions: a
priori estimates on the firing rate, exponential convergence of the
solutions to the steady state when the connectivity parameter is small
enough, and some obstructions to the existence of time-periodic
solutions in order to have insights in their potential domain of
existence in terms of the model parameters.

In conclusion, we complete the mathematical analysis of the NNLIF
system proving global-in-time existence for the delayed NNLIF model
and making progress in the study of its long time behaviour.  In this
way we contribute to a better understanding of the NNLIF model, with a
rich variety of phenomena, able to reproduce biological facts.  In the
light of these and earlier results we can conclude that the blow-up
phenomenon, observed when the synaptic delay is neglected, appears due
to this simplification.  However, when the transmission delay is
included in the model, solutions tend to the unique steady state
when the connectivity parameter $b$ is small enough. For $b$ large
enough so that no steady states exist, we also prove that there
are no periodic solutions and it seems that the firing rate increases
in time, as it was observed numerically in \cite{CS17-2}. The
methods developed in this article could also be used to prove similar
results for the NNLIF model with delay and refractory state for both;
one \cite{CP} and two populations \cite{CS17, CS17-2}.

Several questions regarding the NNLIF model remain open:
what happens with the  solutions of the non-delayed NNLIF model 
after a blow-up phenomenon,
the convergence toward the stationary state for average-inhibitory
case (when entropy methods break), 
the existence and stability of periodic solutions
and the analysis of  possible multistability phenomena.

\appendix

\section{Technical results from the literature}\label{TeXniK}
For completeness we include in this appendix 
two theorems often used in the study of NNLIF models: 
a Poincar\'{e}-like inequality and an entropy equality.
\begin{theorem}[Poincar\'{e}-like inequality]
There exists $\eta>0$ such that for every 
$b\in[-\eta,\eta]$, there exists $\gamma>0$ depending on 
$b$ such that for any measurable function $h$ satisfying 
$\int_{-\infty}^{V_F}\rho_\infty(v)h(v)=1$, the following 
inequality holds:
	$\gamma \int_{-\infty}^{V_F} \rho_\infty(v) (h(v)-1)^2dv \leqslant 
\int_{-\infty}^{V_F}\rho_\infty(v)
\left[ \dfrac{\partial h}{\partial v} \right]^2(v)dv$,
where $\rho_\infty$ is the steady state associated with $b$.
\end{theorem}
The proof of this result was done in \cite{CCP} for $b=0$,
and it is easily extended for small enough $b$, because it only uses 
the behavior of $\rho_\infty$ around $V_F$ and $-\infty$, as
it was explained  in \cite{carrillo2014qualitative}.
That behavior does not change asymptotically when 
$b$ goes from 0 to a small value.

In the rest of the appendix we prove the following entropy equality,
where we assume $b_0=0$, without loss of generality.
\begin{theorem}[Entropy equality]\label{thm:entropy}
For all convex function $G:\R_+\to\R$ of class $C^2$ 
and any initial condition $(\rho^0,N^0)$ satisfying the 
hypotheses of Theorem \ref{thm:globalexist}, if there exists $C_0>0$ such that
$\rho^0(v)\leq C_0 \rho_{\infty}(v)$, for all $v\in (-\infty,V_F]$, 
	then for all $t>0$, there exists $C(t)$ such that for all $v\in(-\infty,V_F],\ \rho(v,t)\leq C(t) \rho_{\infty}(v) $ and 
the corresponding classical solution $\rho$ of \eqref{delay_original}-\eqref{condiciones1}-\eqref{condiciones2} satisfies 	
	\begin{multline}
	\displaystyle{ \dfrac{d}{d t}\int_{-\infty}^{V_F} \rho_\infty 
G\left(\dfrac{\rho}{\rho_\infty}\right) d v } \ 
 = 
	\displaystyle{ -a\int_{-\infty}^{V_F}\rho_\infty 
\left[\dfrac{\partial}{\partial v}\left(\dfrac{\rho}{\rho_\infty}\right)
\right]^2G''\left(\dfrac{\rho}{\rho_\infty}\right)d v   } 
\\ 
	- \displaystyle{N_{\infty}
\Big[G\left(\dfrac{N}{N_\infty}\right) - G\left(\dfrac{\rho}{\rho_\infty}\right) - \Big( \dfrac{N}{N_\infty} - \dfrac{\rho}{\rho_\infty} \Big)
G'\left(\dfrac{\rho}{\rho_\infty}\right)   \Big]\bigg\rvert_{v=V_R}  }
\\
	+  \displaystyle{b(N(t-D)-N_\infty)\int_{-\infty}^{V_F}\dfrac{\partial \rho_\infty}{\partial v}\Big[G\left(\dfrac{\rho}{\rho_\infty}\right) -\dfrac{\rho}{\rho_\infty}G'\left(\frac{\rho}{\rho_\infty}\right)\Big]d v }.
\label{eq: entropy}
	\end{multline}
	\end{theorem}
{\bf Proof.}
According to Section \ref{sec: local}, for $a=1$, $V_F=0$ and $t<D$, if we denote
\[ \mathcal{G}(x,\tau,\xi,\eta) = \dfrac{1}{\sqrt{4\pi (\tau-\eta) }} e^{-\dfrac{(x-\xi)^2}{4(\tau-\eta)}}   \]
the heat kernel, the solution $\rho$ satisfies
\begin{multline}
\rho(v,t)=e^{2t}\int_{-\infty}^{0}\mathcal{G}
\left(e^t v - \int_0^{\frac{1}{2}(e^{2t}-1)} e^{-s}\mu(s-D) ds,\frac{1}{2}(e^{2t}-1),\xi,0\right)\rho^0(e^t \xi )\ d\xi\\ 
- e^t\int_0^{\frac{1}{2}(e^{2t}-1)} M(\eta)\mathcal{G}
\left(e^t v - \int_0^{\frac{1}{2}(e^{2t}-1)} e^{-s}\mu(s-D)d s ,\frac{1}{2}(e^{2t}-1),s(\eta),\eta\right) \ d\eta
\\ 
+ e^t \int_0^{\frac{1}{2}(e^{2t}-1)} M(\eta)\mathcal{G}\left(e^t v - \int_0^{\frac{1}{2}(e^{2t}-1)} e^{-s}\mu(s-D) d s,\frac{1}{2}(e^{2t}-1),s_1(\eta),\eta\right)\ d\eta,
\end{multline}
where $M$ only depends on $N^0$ as long as $t<D$.\\
Due to the forms of $\mathcal{G}$ and $\rho_{\infty}$, the last two terms decrease at least as fast as $\rho_{\infty}$ at $-\infty$. The first term satisfies
\begin{multline}
e^{2t}\int_{-\infty}^{0}\mathcal{G}\Big(e^t v - \int_0^{\frac{1}{2}(e^{2t}-1)} e^{-s}\mu(s-D) ds,\frac{1}{2}(e^{2t}-1),\xi,0\Big)\rho^0(e^t \xi )\ d\xi\\ 
\leq e^{2t}C_0\int_{-\infty}^{0}\mathcal{G}\Big(e^t v - \int_0^{\frac{1}{2}(e^{2t}-1)} e^{-s}\mu(s-D) ds,\frac{1}{2}(e^{2t}-1),\xi,0\Big)\rho_\infty(e^t \xi )\ d\xi, \\
= e^{t}C_0\int_{-\infty}^{0}\mathcal{G}\Big( e^tv - \int_0^{\frac{1}{2}(e^{2t}-1)} e^{-s}\mu(s-D) ds,\frac{1}{2}(e^{2t}-1),e^{-t}\xi,0\Big)\rho_\infty(\xi )\ d\xi.
\end{multline}
Hence, there exists $C(t)$ such that
$\rho(v,t) \leqslant C(t) \rho_{\infty}(v)$,
$\forall v\in(-\infty,V_F]$.  
As we said previously, it is then true for all $a,V_F,t$ using changes of variables and successive time intervals of length less than $D$.

To conclude we prove the entropy equality \ref{eq: entropy}.
First, we calculate $\dfrac{\partial}{\partial v} \left(\dfrac{\rho}{\rho_\infty}\right)$:
\[
\dfrac{\partial}{\partial v} \left(\dfrac{\rho}{\rho_\infty}\right) = \dfrac{1}{\rho_\infty}\dfrac{\partial \rho}{\partial v}-\dfrac{\rho}{\rho_\infty^2}\dfrac{\partial \rho_\infty}{\partial v}
\]
and thus 
\[ \dfrac{\partial \rho}{\partial v} = \rho_\infty\dfrac{\partial}{\partial v} \left(\dfrac{\rho}{\rho_\infty}\right)+\dfrac{\rho}{\rho_\infty}\dfrac{\partial \rho_\infty}{\partial v}
,\]
and  
\[\dfrac{\partial^2 \rho}{\partial v^2} = \rho_\infty\dfrac{\partial^2}{\partial v^2} \left(\dfrac{\rho}{\rho_\infty}\right)+2\dfrac{\partial \rho_\infty}{\partial v} \dfrac{\partial}{\partial v} \left(\dfrac{\rho}{\rho_\infty}\right) + \dfrac{\rho}{\rho_\infty}\dfrac{\partial^2 \rho_\infty}{\partial v^2}.
\] 
These expressions allow us to do the following computations 
\[
\begin{array}{rcl}
\displaystyle{ \dfrac{\partial}{\partial t}\left[ \rho_\infty 
G\left(\dfrac{\rho}{\rho_\infty}\right) \right] }&=&\displaystyle{  \dfrac{\partial \rho}{\partial t}G'\left(\dfrac{\rho}{\rho_\infty}\right)}\\
&=& \displaystyle{  \left(\dfrac{\partial}{\partial v}\Big[\big(v-bN(t-D)\big)\rho\Big]+a\dfrac{\partial^2 \rho}{\partial v^2} + \delta_{V_R}  N\right)G'\left(\dfrac{\rho}{\rho_\infty}\right)}\\
&=& \displaystyle{  \Bigg(v \rho_\infty\dfrac{\partial}{\partial v} \left(\dfrac{\rho}{\rho_\infty}\right)+v\dfrac{\rho}{\rho_\infty}\dfrac{\partial \rho_\infty}{\partial v}+\rho}
\\
&&
\displaystyle{-bN(t-D) \rho_\infty\dfrac{\partial}{\partial v} \left(\dfrac{\rho}{\rho_\infty}\right)-bN(t-D)\dfrac{\rho}{\rho_\infty}\dfrac{\partial \rho_\infty}{\partial v}  + \delta_{V_R}  N}\\
&& \displaystyle{+a\rho_\infty\dfrac{\partial^2}{\partial v^2} \left(\dfrac{\rho}{\rho_\infty}\right) +2a\dfrac{\partial \rho_\infty}{\partial v} \dfrac{\partial}{\partial v} \left(\dfrac{\rho}{\rho_\infty}\right) + a\dfrac{\rho}{\rho_\infty}\dfrac{\partial^2 \rho_\infty}{\partial v^2} \Bigg)G'\left(\dfrac{\rho}{\rho_\infty}\right)}\\
&=& \displaystyle{    \Bigg(  \left(v\rho_\infty+2a\dfrac{\partial \rho_\infty}{\partial v}\right) \dfrac{\partial}{\partial v} \left(\dfrac{\rho}{\rho_\infty}\right)+a\rho_\infty\dfrac{\partial^2}{\partial v^2} \left(\dfrac{\rho}{\rho_\infty}\right)}\\
&&\displaystyle{  +\dfrac{\rho}{\rho_\infty}\Big[ v
\dfrac{\partial \rho_\infty}{\partial v} + \rho_\infty  + a\dfrac{\partial^2 \rho_\infty}{\partial v^2}  \Big]  -bN(t-D)\dfrac{\rho}{\rho_\infty}\dfrac{\partial \rho_\infty}{\partial v}     } \\
&&\displaystyle{    -bN(t-D) \rho_\infty\dfrac{\partial}{\partial v} \left(\dfrac{\rho}{\rho_\infty}\right)  + \delta_{V_R}  N \Bigg)G'\left(\dfrac{\rho}{\rho_\infty}\right)}.
\end{array}
\]
As $\rho_\infty$ is a steady state, we have
\begin{multline}
\dfrac{\rho}{\rho_\infty}\Big[ v
\dfrac{\partial \rho_\infty}{\partial v} + \rho_\infty  + a\dfrac{\partial^2 \rho_\infty}{\partial v^2}  \Big]  -bN(t-D)\dfrac{\rho}{\rho_\infty}\dfrac{\partial \rho_\infty}{\partial v}\ 
\\
=
\dfrac{\rho}{\rho_\infty}\Big[ bN_\infty\dfrac{\partial \rho_\infty}{\partial v} -\delta_{V_R}N_{\infty}  \Big]  -bN(t-D)\dfrac{\rho}{\rho_\infty}\dfrac{\partial \rho_\infty}{\partial v}\\
=\ b(N_\infty-N(t-D))\dfrac{\rho}{\rho_\infty}\dfrac{\partial \rho_\infty}{\partial v} -\delta_{V_R}N_{\infty}\dfrac{\rho}{\rho_\infty},
\end{multline}
and combining the two previous computations 
\[
\begin{array}{rcl}
\displaystyle{ \dfrac{\partial}{\partial t}\left[ \rho_\infty G\left(\dfrac{\rho}{\rho_\infty}\right) \right] } &=& \displaystyle{    \Bigg(  \left(v\rho_\infty+2a\dfrac{\partial \rho_\infty}{\partial v}\right) \dfrac{\partial}{\partial v} \left(\dfrac{\rho}{\rho_\infty}\right)+a\rho_\infty\dfrac{\partial^2}{\partial v^2} \left(\dfrac{\rho}{\rho_\infty}\right)}\\
&&\displaystyle{  +\delta_{V_R}N_{\infty}\Big( \dfrac{N}{N_\infty} - \dfrac{\rho}{\rho_\infty} \Big) + b(N_\infty-N(t-D))\dfrac{\rho}{\rho_\infty}\dfrac{\partial \rho_\infty}{\partial v}  } \\
&&\displaystyle{    -bN(t-D) \rho_\infty\dfrac{\partial}{\partial v} \left(\dfrac{\rho}{\rho_\infty}\right) \Bigg)G'\left(\dfrac{\rho}{\rho_\infty}\right)}.
\end{array}
\]
On the other hand, since 
$\dfrac{\partial}{\partial v}G\left(\frac{\rho}{\rho_\infty}\right) = 
\dfrac{\partial}{\partial v}\left(\dfrac{\rho}{\rho_\infty}\right)
G'\left(\frac{\rho}{\rho_\infty}\right) 
$
and
\[ \dfrac{\partial^2}{\partial v^2}G\left(\frac{\rho}{\rho_\infty}\right) = 
\left[\dfrac{\partial}{\partial v}\left(\dfrac{\rho}{\rho_\infty}\right)\right]^2G''\left(\dfrac{\rho}{\rho_\infty}\right) +  \dfrac{\partial^2}{\partial v^2}\left(\dfrac{\rho}{\rho_\infty}\right)G'\left(\frac{\rho}{\rho_\infty}\right),\]
we obtain
\[
\begin{array}{rcl}
\displaystyle{ \dfrac{\partial}{\partial t}\left[ \rho_\infty G\left(\dfrac{\rho}{\rho_\infty}\right) \right] } &=& \displaystyle{   \left(v\rho_\infty+2a\dfrac{\partial \rho_\infty}{\partial v}\right) \dfrac{\partial}{\partial v}G\left(\frac{\rho}{\rho_\infty}\right)+a\rho_\infty\dfrac{\partial^2}{\partial v^2} 
G\left(\frac{\rho}{\rho_\infty}\right)}\\
&& \displaystyle{ -a\rho_\infty \left[\dfrac{\partial}{\partial v}\left(\dfrac{\rho}{\rho_\infty}\right)\right]^2G''\left(\dfrac{\rho}{\rho_\infty}\right)  +\delta_{V_R}N_{\infty}\Big( \dfrac{N}{N_\infty} - \dfrac{\rho}{\rho_\infty} \Big)
G'\left(\dfrac{\rho}{\rho_\infty}\right)   } \\
&&\displaystyle{  + b(N_\infty-N(t-D))\dfrac{\partial \rho_\infty}{\partial v}\dfrac{\rho}{\rho_\infty}G'\left(\frac{\rho}{\rho_\infty}\right)  -bN(t-D) \rho_\infty\dfrac{\partial}{\partial v} G\left(\frac{\rho}{\rho_\infty}\right) }.
\end{array}
\]
Then, using the equation for $\rho_\infty$, we write
\[
\begin{array}{rcl}
\displaystyle{ \dfrac{\partial}{\partial t}\left[ \rho_\infty G\left(\dfrac{\rho}{\rho_\infty}\right) \right] } &=& \displaystyle{ \dfrac{\partial}{\partial v}\Big( v\rho_\infty G\left(\frac{\rho}{\rho_\infty}\right)\Big) + a\dfrac{\partial^2}{\partial v^2}\Big(\rho_\infty G\left(\frac{\rho}{\rho_\infty}\right)\Big)  }\\
&& \displaystyle{  +\delta_{V_R}N_{\infty}\Big[ \Big( \dfrac{N}{N_\infty} - \dfrac{\rho}{\rho_\infty} \Big)G'\left(\dfrac{\rho}{\rho_\infty}\right) + 
G\left(\dfrac{\rho}{\rho_\infty}\right)  \Big]  }\\
&&\displaystyle{ -a\rho_\infty \left[\dfrac{\partial}{\partial v}\left(\dfrac{\rho}{\rho_\infty}\right)\right]^2G''\left(\dfrac{\rho}{\rho_\infty}\right) -bN(t-D) \rho_\infty\dfrac{\partial}{\partial v} G\left(\frac{\rho}{\rho_\infty}\right)   } \\
&&\displaystyle{   -bN_\infty\dfrac{\partial \rho_\infty}{\partial v}
G\left(\dfrac{\rho}{\rho_\infty}\right)  + b(N_\infty-N(t-D))\dfrac{\partial \rho_\infty}{\partial v}\dfrac{\rho}{\rho_\infty}G'\left(\frac{\rho}{\rho_\infty}\right)   }.
\end{array}
\]
Finally, as $\rho(V_F,t)=\rho_\infty(V_F)=0$ and $\frac{\partial \rho_\infty}{\partial v}(V_F)=-aN_\infty\neq 0$, we can apply l'H\^{o}pital's rule :
\[
\forall t\in\R_+, \ \ \ \lim_{v\to V_F}\dfrac{p(v,t)}{\rho_\infty(v)}=\lim_{v\to V_F}\dfrac{\frac{\partial \rho}{\partial v}(v,t)}{\frac{\partial \rho_\infty}{\partial v}(v)} = \dfrac{N(t)}{N_\infty}, 
\]
and integrating 
we obtain
$\int_{-\infty}^{V_F}  \dfrac{\partial}{\partial v}\Big( v\rho_\infty G\left(\frac{\rho}{\rho_\infty}\right)\Big) d v = \Big[ v\rho_\infty G\left(\frac{\rho}{\rho_\infty}\right)\Big]_{-\infty}^{V_F}=0$
and
\[
\begin{array}{rcl}
\displaystyle{\int_{-\infty}^{V_F}a\dfrac{\partial^2}{\partial v^2}\Big(\rho_\infty G\left(\frac{\rho}{\rho_\infty}\right)\Big)d v }&=&\displaystyle{ \Big[a\dfrac{\partial}{\partial v}\left(\rho_\infty G\left(\frac{\rho}{\rho_\infty}\right)\right)\Big]_{-\infty}^{V_F}}\\
&=&\displaystyle{\left[a\left(\dfrac{\partial \rho}{\partial v} - 
\dfrac{\rho}{\rho_\infty}\dfrac{\partial \rho_\infty}{\partial v} \right) G'\left(\frac{\rho}{\rho_\infty}\right) + 
a\dfrac{\partial \rho_\infty}{\partial v}G\left(\frac{\rho}{\rho_\infty}\right)
\right]_{-\infty}^{V_F}}\\
&=&\displaystyle{ \Big(-N + \dfrac{N}{N_\infty}N_\infty \Big) G'\left(\frac{N}{N_{\infty}}\right) - N_\infty G\left(\dfrac{N}{N_\infty}\right)  }\\
&=& - N_\infty G\left(\dfrac{N}{N_\infty}\right), 
\end{array}
\]
where every integral is defined and finite thanks to the inequality $\rho(v,t)\leq C(t) \rho_{\infty}(v)$.\\
Eventually, we have, by integration by parts and using boundary conditions,
\[ \int_{-\infty}^{V_F}bN(t-D) \rho_\infty\dfrac{\partial}{\partial v} G\left(\frac{\rho}{\rho_\infty}\right)d v = -bN(t-D) \int_{-\infty}^{V_F} \dfrac{\partial \rho_\infty}{\partial v} G\left(\frac{\rho}{\rho_\infty}\right)d v \]
and the result comes from this integral and the previous computations.
\qed

\bibliographystyle{siam}
\bibliography{Biblio}
\end{document}